\documentclass[12pt]{amsart}
\usepackage{a4wide,graphicx}
\usepackage{color,enumerate}
\usepackage{soul}
\usepackage{euscript}

 \usepackage{caption}
 \usepackage[labelformat=simple,labelfont={}]{subcaption}

\let\eps\varepsilon  
\newcommand{\Nb}{{\mathbb N}}  
\newcommand{\Rb}{{\mathbb R}} 
\newcommand{\diver}{\operatorname{div}}  
\newcommand{\supp}{\operatorname{supp}}  
\newcommand{\sspan}{\operatorname{span}} 
\newcommand{\di}{\displaystyle}

\newcommand{\TM}{{\EuScript T}_M}

\newcommand{\ub}{\mathbf{u}}
\newcommand{\nb}{\mathbf{0}}

\newtheorem{theorem}{Theorem} 
\newtheorem{lemma}{Lemma}    
\newtheorem{proposition}{Proposition}   
\newtheorem{remark}{Remark}    
\newtheorem{corollary}{Corollary}   
\newtheorem{definition}{Definition}

\begin{document}  

\title[Unsaturated flow in porous media with dynamic capillary pressure]
{The unsaturated flow in porous media with dynamic capillary pressure}

\author[J. P. Mili\v{s}i\'{c}]{Josipa-Pina Mili\v{s}i\'{c}}
\address{University of Zagreb, Faculty of Electrical Engineering and Computing, Unska 3,
10000 Zagreb, Croatia}
\email{pina.milisic@fer.hr}

\date{\today}


\begin{abstract}
 In this paper we consider a degenerate pseudoparabolic equation
 for the wetting saturation of an unsaturated two-phase flow in porous media with dynamic capillary pressure-saturation
 relationship where the relaxation parameter depends on the saturation. Following the approach given in \cite{M10}
 the existence of a weak solution is proved using Galerkin approximation and regularization techniques. 
 A priori estimates needed for passing to the limit when the regularization parameter goes to zero
 are obtained by using appropriate test-functions, motivated by the fact that considered PDE allows a natural 
 generalization of the classical Kullback entropy. Finally, a special care was given in obtaining an estimate of the mixed-derivative
 term by combining the information from the capillary pressure with the obtained a priori estimates on the saturation. 
\end{abstract}

  

\maketitle


\section{Introduction}\label{Intro}

 Two-phase flow processes in porous media appear in many real situations, 
 such as unsaturated water flow in the subsurface, remediation of contaminated sites and water-oil displacement during oil recovery.
 The physical-mathematical model of the two-phase flow 
 consists of two mass balance equations, Darcy relations for the two phases and the capillary condition 
  which is a constitutive relationship between the wetting-phase saturation and capillary pressure $P_c$.
 In this paper we include the dynamic 
 capillary pressure relation proposed by Hassanizadeh and Grey \cite{HG93}
 to the standard two-phase flow equations. This results in a transport equation containing higher-order mixed derivatives.
 Moreover, the obtained model represents degenerate pseudoparabolic equation
 and our goal is to prove the global-in-time existence of its weak solutions.
 More specifically, we consider the equation
 \begin{align}
  \partial_t S = \diver \Big( a(S) \nabla \big( \partial_t \beta(S) - P_c(S) \big) \Big),
 \label{Mil.0}
\end{align}
where
\begin{align*}
  a(0) & = a(1) = 0,\\
  -P_c'(S) \to  +\infty \; & \textrm{ when  }\; S \to 0, \; \textrm{ and } \; S \to 1,
\end{align*}
\begin{align*}
  \beta(S) & = \int_0^S \tau(s) ds,\\
 \tau(0)  = 0,\; \tau(S) \textrm{ is bounded when } & S \to 1 \textrm{ and } 
 \tau \textrm{ is monotone increasing. }
\end{align*}
The solution of equation \eqref{Mil.0} satisfies certain initial-boundary conditions which will be specified later.

 We note that the nonequilibrium capillary effects in problems of enhancing oil and gas recovery from
rocks were proposed in classical book by Barenblatt, Entov and Ryzhik \cite{BER90}, and then investigated by 
numerous scientists (see for example \cite{B71, BV80, BPS02, V78}, \cite{BeHas01, Hass97}).
The mathematical analysis of dynamic capillary pressure models for Richards' equation was intensively studied.
A significant contribution has been done by Mikeli\'c \cite{M10} where the existence, for any time interval,
 of an appropriate weak solution of the Richards' equation 
 was proved.
 The existence of a weak solution to the Richards' equation with dynamic capillary pressure and hysteresis
 is studied by Schweizer \cite{Sch12-2}. For the two-phase flow model  with dynamic capillary pressure and hysteresis,
 the existence of a weak solution is proved by Koch, R\"atz and Schweizer \cite{KoRaSch13} under assumption of non-degenerated mobilities.
  
  The first existence result for the two-phase flow model with dynamic capillary pressure and saturation dependent relaxation parameter
 is obtained by Cao and Pop in \cite{CaoPop15}. They considered the system in global pressure formulation with
 a strictly positive relaxation parameter 
 and a constant Dirichlet boundary condition for the saturation. 
 
  In this work we consider the two-phase flow model with dynamic capillary pressure 
 and degenerate saturation-dependent relaxation parameter. Moreover, we impose mixed boundary
 conditions with a non-constant boundary saturation.
 As our result shows (see also \cite{M10} and \cite{CaoPop15}), the existence theorem can be proved under certain
 relations between the orders of zeros of relative permeabilities
  and relaxation parameter and the order of singularities of capillary pressure function (see Proposition~3).
  Concerning these relations it is sufficient to analyze
  the  case of the countercurrent imbibition flow instead of the full two-phase flow system.
  Therefore, we consider the countercurrent imbibition flow in the two-phase system that changes wettability (see \cite{DonAl08}).
 
\vspace{2mm}
 
{\bf Key ideas.} 
 The motivation for this article comes from the result of Mikeli\'c \cite{M10} on a global existence for the following nonlinear degenerate
 pseudoparabolic equation
 \begin{align}
   \partial_t S  = \diver \Big( k(S) \nabla \big( \tau \partial_t S - P_c(S) \big) \Big).
  \label{Mik.0}
 \end{align}
 Mikeli\'c observed that equation \eqref{Mik.0} allows a natural generalization of the classical Kullback entropy given by
$ \mathcal{E}_M''(s) = 1/k(s)$. Namely, formal calculation gives
\begin{align}
 \frac{d}{dt} \int_\Omega \Big( \mathcal{E}_M(S) + \frac{\tau}{2} |\nabla S|^2  \Big) dx = \int_\Omega P_c'(S)|\nabla S|^2 dx + \textrm{ low order terms. }
 \label{Mik.0E}
\end{align}
 The main difficulties for the mathematical analysis of \eqref{Mik.0} were the degenerations of the coefficients, i.e. $k(0) = 0$ and
$-P_c'(S) \to +\infty$ when $S \to 0$. These degenerate coefficients and presence of the initial and the boundary conditions 
lead to unbounded non-integrable $\mathcal{E}_M'$. Because of that reason, in \cite{M10} the degenerate
coefficients were firstly regularized, by introducing the regularization parameter $\eps$. After that, the entropy estimates necessary for passing to the limit in the
weak formulation when $\eps \to 0$ were derived. 
The crucial step in mentioned existence result was the use of $\mathcal{E}_M'$ as a test-function
for obtaining needed a priori estimates. 

 In comparison with equation \eqref{Mik.0} considered in \cite{M10}, besides the additional degeneration of the coefficient $a(S)$ in $S=1$, 
 and the fact that $-P_c'(S) \to +\infty$ for $S  \to 0$ and $S  \to 1$, the main novelty in this article is the dependence of the relaxation parameter 
 $\tau$ upon the variable $S$,
 i.e. after rewriting, the equation of our interest \eqref{Mil.0} reads
 \begin{align}
   \partial_t S = \diver \Big( a(S) \nabla \big( \tau(S) \partial_t S  - P_c(S) \big) \Big).
  \label{Mil.1}
 \end{align}
 Similarly like in \cite{M10}, we observe that equation \eqref{Mil.1} allows a natural generalization of the classical Kullback
 entropy. For the equation of our interest \eqref{Mil.0}, the entropy functional is given by formula
\begin{align}
 \mathcal{E}(S) = S \int_{S_D}^S \frac{\tau(\xi)}{a(\xi)} d\xi - \int_{S_D}^S \xi \frac{\tau(\xi)}{a(\xi)} d\xi,
 \label{Czisar-Kullback}
\end{align}
where $S_D$ comes from the boundary condition.
Easy calculation gives $\mathcal{E}''(S) = \tau(S)/a(S)$.
 In order to obtain needed a priori estimates, adapting the idea given in \cite{M10}, we use 
\[ \mathcal{E}'(S) =  \int_{S_D}^S \frac{\tau(\xi)}{a(\xi)}  d\xi, \]
as a test-function in a weak formulation given later by \eqref{WF}, leading to
an estimate similar to \eqref{Mik.0E}. 

\section{Physical model and basic equations}

We consider two-phase water-oil flow (wetting and non-wetting immiscible fluids)
in an isotropic and homogeneous porous medium. Following \cite{Bear68}, the generalized Darcy law gives:
 \begin{align}
   \ub_w = -\lambda_w(S) k \nabla p_w,\quad  \ub_o = -\lambda_o(S) k \nabla p_o.
   \label{Darcy}
 \end{align}
Here the subscripts $w$ and $o$ correspond, respectively, to the water (wetting) and the oil (non-wetting) fluids,
$\ub_i$ are the fluxes of the phases, $p_i$ are their pressures, and $\lambda_i$ are phase mobilities, depending on the 
wetting-phase saturation $S = S_w$.
Furthermore, $k$ is the absolute permeability of the porous medium,
and the gravity effects are neglected for simplicity.
The mass conservation laws for both phases have the form:
\begin{align}
 \phi \frac{\partial S}{\partial t} + \diver \ub_w  = 0, \quad
 \phi \frac{\partial(1-S)}{\partial t} + \diver \ub_o  = 0, \label{eq-1-2}
\end{align}
where $\phi$ is the porosity of the medium.
The model \eqref{Darcy}--\eqref{eq-1-2} has to be completed with the capillary pressure law 
which has the form
\[ p_o - p_w = p_c, \]
where, due to \cite{HG93}, the capillary pressure saturation relationship is given by
\begin{align}
 p_c = P_c(S) - \tau(S) \frac{\partial S}{\partial t}. 
 \label{HassGrey}
\end{align}
Here, $P_c(S)$ is static capillary pressure function and $\tau(S)$ is the relaxation parameter.
Expression \eqref{HassGrey}
suggests that a new equilibrium between the dynamic and the static capillary pressure is not attained instantaneously
 and the relaxation parameter $\tau$ depends on the water saturation. 
 Following the work of Cuesta, van Duijn and Hulshof \cite{CueDujHul00} and the references therein \cite{Hass97, SmVacVau71},
 we take that $\tau(S)$ vanishes as $S \to 0$. This assumption can be explained as in \cite{BPS02}. Namely, 
 at low water saturation the water flows to the narrowest pores and the time needed to reconfigure the water distribution tends to zero.
 We note that this degeneration introduces additional mathematical difficulty into the system in comparison
 with \cite{M10}. For purposes of our mathematical analysis (see inequality \eqref{Ineq.1}), 
 we take that $\tau(S)$ is monotone increasing function. 

 Hereafter, we normalize the porosity $\phi$ and the permeability $k$ to one. 
 By adding equations \eqref{eq-1-2} we get 
 \begin{align}
 \frac{\partial S}{\partial t} +  \diver \big( b(S) \ub \big) = \diver\big( a(S) \nabla (-p_c) \big),
\label{eq-5}
\end{align}
 where the coefficients are given by 
\begin{align*}
 b(S) = \frac{\lambda_w(S)}{\lambda_w(S) + \lambda_o(S)},\quad a(S) =   \frac{\lambda_w(S) \lambda_o(S)}{\lambda_w(S) + \lambda_o(S)}.                    
\end{align*}

For the simplicity reasons we consider here the capillary countercurrent imbibition flow (see \cite{BPS02}) where, due
 to the incompressibility of both fluids, the bulk fluid flow is identically equal zero:
 \[ \ub = \ub_o + \ub_w = \nb. \]
More precisely, we study the following equation for the water saturation:
\begin{align*}
 \partial_t S + \diver\Big( a(S) P_c'(S)\nabla S\Big) = \diver \Big( a(S) \nabla \big( \tau(S) \partial_t S  \big) \Big).
\end{align*}
 which is actually the equation of our interest, written in \eqref{Mil.0}.

\section{Equation and its weak formulation}

 Let $\Omega \subset \Rb^n$ be an open, bounded  
 domain, with Lipschitz boundary $\partial \Omega$. The boundary $\partial \Omega$ is decomposed into the Dirichlet
 part $\partial_D \Omega$ and the Neumann part $\partial_N \Omega$, where
 $\partial \Omega = \overline{\partial_D \Omega} \cup \overline{\partial_N \Omega}$ and 
 $\partial_D \Omega \cap \partial_N \Omega = \emptyset$.
 We suppose that the Dirichlet boundary $\partial_D \Omega$ is closed relative to $\partial \Omega$ and has a
 positive $(n-1)$-dimensional Hausdorff measure. The Neumann boundary $\partial_N \Omega$ may be empty.
 Let $Q_T = \Omega \times (0,T)$, $T> 0$. Moreover, we use the notation
 $\Gamma_D = \partial_D \Omega \times (0,T)$ and $\Gamma_N = \partial_N\Omega \times (0,T)$.
 Here $(0,T)$, $T>0$, is the time interval.
  We are looking for a solution to equation \eqref{Mil.0}
which satisfies the following boundary and initial conditions:
  \begin{align}
    S & = S_D \quad  \textrm{ on } \; \Gamma_D, \label{eq-10}\\ 
   a(S) \nabla \big( P_c(S) - \partial_t \beta(S) \big)  \cdot \nu & = \mathcal{R}  \quad  \textrm{ on } \; \Gamma_N, \label{eq-11}\\ 
    S(x, 0) & = S_i(x)  \quad  \textrm{ on } \; \Omega, \label{eq-12}
  \end{align}
  where $\nu$ is the unit outer normal on $\partial\Omega$, $\mathcal{R}$ is a given flux, 
  $S_D$ is the given saturation at $\Gamma_D$ and
 $S_i$ is the initial saturation of the wetting phase.
 We impose following hypotheses on the coefficients and on the data, for $ s \in [0,1]$:
 
 \vskip 2mm
 
 {\bf (H1)} There are constants $\mu > 0$ and $\lambda > 0$
 such that 
\[  a(s) = \frac{s^{\mu} \cdot (1-s)^{\lambda}}{s^\mu + (1-s)^\lambda}, \]
giving $a(0) = a(1) = 0$.

 \vskip 2mm
 
 {\bf (H2)} For the relaxation parameter we suppose that there exist constant  $\gamma > 0$ and $\TM > 1$,
 such that
 \begin{align}
              \tau(s) = \frac{s^{\mu}}{s^\mu + (1-s)^\lambda} 
              \Big[  \TM  + \frac{(1-s)^{\lambda}}{s^\gamma} \Big].
              \label{Def.tau.1}
            \end{align}
Furthermore, we take $\mu > \gamma \geq 0$ so we get 
$ \tau(0) = 0$ and $\tau(1) =  \TM$. The condition $\TM > 1$ together with $\mu > \gamma$ ensures that the relaxation parameter $\tau$
is increasing function of water saturation $s$.

\vskip 2mm

 {\bf (H3)} For the derivative of the capillary pressure we suppose that there exist constants
 \[  0<\beta_1 \leq \gamma < \mu\; \textrm{ and }\; 0<\beta_2\leq \lambda \]
 such that
 \[ -P_c'(s) =  \frac{g(s)}{s^{\beta_1}} + \frac{h(s)}{(1-s)^{\beta_2}}, \]
 where $g, h \in C^\infty_c(\Rb)$ are arbitrary functions strictly positive on $[0,1]$.
 Note that 
 \[ -P_c'(s) \to +\infty \; \textrm{ when } \; s \to 0^+,\; s \to 1^-. \] 
 
 \vskip 2mm
 
  {\bf (H4)} We assume that $S_D \in C^1([0,T];H^1(\Omega))$, $0 < S_{D_{min}} \leq S_D(x,t) \leq S_{D_{max}} < 1$ a.e. on $Q_T$.
  Moreover, we assume that $\di \|\partial_t S_D\|_{L^\infty(Q_T)} < +\infty$.
  The initial condition $S_i$ belongs to $H^1(\Omega)$ and $0 \leq S_i \leq 1$ a.e. on $\Omega$.
 
  \vskip 2mm
 
  {\bf (H5)} We assume that the boundary-flux $\mathcal{R}$ satisfies
  \begin{align}
   \mathcal{R}(x,t,s) = R_0(x,t) \sigma(s), \; \textrm{ where } \; R_0 \in C^1(\overline{\Gamma}_N),\;
   \sigma \in C_c^\infty(0,1).
   \label{BoundFlux}
  \end{align}

  \vskip 2mm
   
   {\bf (H6)} The initial saturation satisfies the finite entropy condition:
  \begin{align}
   \gamma > 2,\; \lambda > 2 \; \textrm{ and } \; \int_\Omega S_i^{2-\gamma}(x) dx + \int_\Omega\big( 1-S_i(x)  \big)^{2-\lambda} dx < +\infty.
   \label{FEC}
  \end{align}

   \vskip 2mm
  
 \begin{remark}
  In order to avoid technical difficulties we use explicit forms for $a$ and $\tau$. This can be easily generalized,
 for example using the strictly positive on $[0,1]$ function $f \in C^\infty_c(\Rb)$
  such that
 \[ a(s) = \frac{s^{\mu} \cdot (1-s)^{\lambda}}{f(s)\big( s^\mu + (1-s)^\lambda \big)},\; \textrm{  and } \; 
  \tau(s) = \frac{s^{\mu}}{f(s) \big( s^\mu + (1-s)^\lambda \big)} 
              \Big[  \TM  + \frac{(1-s)^{\lambda}}{s^\gamma} \Big], \]
 where $f$ should be chosen in a way that the relaxation parameter $\tau(s)$ is monotone increasing function.             
 \end{remark}
 The function $\beta(S)$ is defined as 
\begin{align}
 \beta(S) = \int_0^S \tau(s)\, ds,
 \label{beta}
\end{align}
and it is bounded: $0\leq \beta(S) \leq \beta(1)$.
We note that we extend $\tau$ to $(-\infty,0)$ and to $(1,+\infty)$ as follows:
\begin{align*}
 \tau(s) = \begin{cases}
            0, & s< 0,\\
            \tau(1), & s > 1.
           \end{cases}
\end{align*}
So we have
\begin{align*}
 \beta(s) = \begin{cases}
            0, & s< 0,\\
            \beta(1) + \tau(1)(s-1), & s > 1.
           \end{cases}
\end{align*}

From hypotheses  {\bf (H1)}--{\bf (H3)} it is easy to see that the following coefficients are bounded:
\begin{align}
	 \Big\| \frac{a P'_{c}}{\tau} \Big\|_\infty < +\infty, \quad 
	\|a P_c'\|_\infty <  +\infty, \quad 
        \Big\| \frac{a}{\tau} \Big\|_\infty < +\infty, \label{eq-coeff-bounds}
\end{align}
\begin{align}
	 \Big\| \frac{a \tau'}{\tau} \Big\|_\infty < +\infty, \quad 
	\|a \tau^2\|_\infty <  +\infty, \quad 
        \Big\| \frac{\sqrt{a} \tau'}{\tau} \Big\|_\infty < +\infty. \label{eq-coeff-bounds_1}
\end{align}

 Our goal is to obtain a global existence of a weak solution for problem \eqref{Mil.0}--\eqref{eq-10}--\eqref{eq-11}--\eqref{eq-12}, for any time interval.
 Firstly, we define the weak formulation of the problem.
 \begin{definition}\label{WF_Def}
  Let
  \[ V = \{ S \in H^1(\Omega) \colon S_{\vert_{\partial_D \Omega}} = 0 \}. \]
 \end{definition}
  Find $ S \in L^2(Q_T)$, $0 \leq S(x,t) \leq 1$ a.e. on  $Q_T$ such that $\beta(S) \in H^1(Q_T)$,
  $\beta(S)-\beta(S_D) \in L^2(0,T;V)$, $a(S)\nabla \partial_t \beta(S) \in L^2(Q_T)$
  and that
  for all $\varphi$ in $H^1(0,T;V)$, $\varphi_{\vert_{t=T}} = 0$ it holds
  \begin{equation}
   \begin{split}
     -\int_0^T\int_{\Omega} & S \frac{\partial \varphi}{\partial t} dxdt- \int_{\Omega} S_i(x) \varphi(x,0)dx +
   \int_0^T\int_{\partial_N \Omega} \mathcal{R} \varphi d\Gamma dt \\ 
     & - \int_0^T \int_{\Omega} a(S)P_c'(S)\nabla S\cdot \nabla \varphi dx dt 
      + \int_0^T \int_{\Omega} a(S) \nabla \partial_t \beta(S) \cdot \nabla \varphi dx dt = 0.
   \end{split}
   \label{WF}
 \end{equation}
 
 In the following theorem we state the main result of the paper.

\begin{theorem} {\sf Global existence of a weak solution.}\\
 Let $n \leq 3$ and let us suppose hypotheses {\bf (H1)}--{\bf (H6)}.
 Under the assumptions on the exponents $\beta_1$, $\gamma$, $\mu$, and $\lambda$ given by
 Proposition~\ref{Prop_ThirdOrder}, there exists a weak solution $S$ for problem
 \eqref{Mil.0}-\eqref{eq-10}-\eqref{eq-11}-\eqref{eq-12} such that
 $0 \leq S \leq 1$.
 \label{Thm-Main}
\end{theorem}
 
 The rest of the paper is organized as follows. In Section~\ref{Sect.4} we firstly regularize the coefficients
  and define a weak formulation of the regularized problem. Then by Theorem~\ref{Thm-Reg} we state 
  the global existence of a weak solution for the regularized problem.
  Section~\ref{Sect.5} gives proof of Theorem~\ref{Thm-Reg}. More precisely, we start with the Galerkin method for the regularized problem, 
 then in Theorem~\ref{AprEstThm}, Section~\ref{Sect.6} we give a priori estimates uniform with respect to the 
 regularization parameter $\eps$. Special care was given in obtaining a priori estimate of the mixed-derivative term
 in Section~\ref{Sect.7}, where we combined informations from the capillary pressure with those obtained from the entropy functional. 
 Finally, in Section~\ref{Proof_Thm_Main} we pass to the limit when $\eps \to 0$ finishing the proof of Theorem~\ref{Thm-Main}.
 
  \section{Regularization of the coefficients} \label{Sect.4}
 For proving existence of a weak solution given by Definition~\ref{WF_Def}
we need to define a regularized problem. For that purpose
we regularize the coefficients $a$, $P_C'$ and $\tau$ and extend them to the whole $\Rb$ as follows. Let us introduce 
the cutting-function $Z$ defined for $s\in \mathbb{R}$  as 
\begin{align}
  Z(s) = \left\{
\begin{array}{rl}
0 & \text{if } s < 0,\\
s & \text{if } 0 \leq s \leq 1,\\
1 & \text{if } s > 1,
\end{array} \right.
 \label{Cutt.1}
\end{align}
and, for $\eps > 0$, the function 
 \[ Z_\eps(s) = (1-2\eps)Z(s) + \eps. \]
Then we define 
 \begin{equation}
  \begin{split}
   a_{\eps}(s) = a\big( Z_\eps(s)\big),\quad  P_{c,\eps}'(s) = P_c'\big( Z_\eps(s) \big),\quad
   \tau_\eps(s) = \tau \big( Z_\eps(s) \big),
   \end{split}
   \label{eq-13}
  \end{equation}
  where $s\in \Rb$.
We also define 
 \[ \beta_\eps(S) = \int_0^S \tau_\eps(\xi) d\xi. \]
 It is easy to see that the regularized coefficients are bounded, i.e.
 \begin{align}
  0 < m_a^\eps \leq a_\eps(s) \leq \|a\|_\infty,\quad 0 < m_{aP}^\eps \leq a_\eps(s) P_{c_{\eps}}'(s) \leq \|a P_c'\|_\infty, \label{eq-14} \\
  0 < m_\tau^\eps \leq   \tau_\eps(s) \leq \|\tau\|_\infty, \quad 
  m_\tau^\eps s\leq \beta_\eps(s)\leq \EuScript{T}_M s. \label{eq-15}
 \end{align}
We note that constants $m_a^\eps$, $m_{aP}^\eps$, $m_\tau^\eps$ depend on $\eps$. 
Moreover, $m_\tau^\eps = \tau(\eps)$.
Furthermore, we have 
 \begin{align}
  \left| \frac{a_\eps P'_{c_{\eps}}}{\tau_\eps} (s)\right| & \leq \Big\| \frac{a P'_{c}}{\tau} \Big\|_\infty < +\infty, \label{eq-14.1}\\
 0\leq  \frac{a_\eps}{\tau_\eps} (s) & \leq \Big\| \frac{a}{\tau} \Big\|_\infty < +\infty. \label{eq-15.1}
 \end{align}

  Next, we introduce the notion of a weak solution for the regularized version of problem
 \eqref{Mil.0}-\eqref{eq-10}-\eqref{eq-11}-\eqref{eq-12}.
\begin{definition}
 Find $ S_\eps \in H^1(Q_T)$ such that $S_\eps-S_D \in L^2(0,T;V)$, $ \nabla \partial_t \beta (S_\eps) \in L^2(Q_T)$
 and satisfying 
  \begin{equation}
   \begin{split}
    & \int_0^T\int_{\Omega}  (\partial_t S_\eps) \varphi dxdt 
   \int_0^T\int_{\partial_N \Omega} \mathcal{R} \varphi d\Gamma dt \\ 
     & - \int_0^T \int_{\Omega} a_\eps(S_\eps)P_{c_{\eps}}'(S_\eps)\nabla S_{\eps} \cdot \nabla \varphi dx dt 
     + \int_0^T \int_{\Omega} a_\eps(S_\eps) \nabla \partial_t \beta_\eps (S_\eps) \cdot \nabla \varphi dx dt = 0,
  \label{WFReg}
    \end{split}
    \end{equation}
 for all $\varphi \in H^1(0,T;V)$ such that $\varphi_{\vert_{t=T}} = 0$, and $S_\eps(0) = S_i$ in $L^2(\Omega)$.  
\end{definition}

\begin{theorem} {\sf Global existence of a weak solution for the regularized problem.}\\
  Under hypotheses {\bf (H1)}--{\bf (H5)} there exists a solution for the regularized problem
 \eqref{WFReg}.
 \label{Thm-Reg}
\end{theorem}
\section{Proof of Theorem \ref{Thm-Reg}}\label{Sect.5}

 {\sf Step 1: Galerkin approximation.} 
 
 Let $(e_j)_{j\in \Nb}$ be a basis of $V$ and let $V_N = \sspan \{e_1,\ldots,e_N\}$. 
 The Galerkin approximation for the regularized problem \eqref{WFReg}
 reads as follows:
 
 find $\di \beta_\eps(S_{\eps,N}) = \beta_\eps(S_D) + \sum_{j=1}^N \gamma_j(t) e_j(x)$, $\gamma_j \in C^1([0,T])$ such that
 \begin{equation}
  \begin{split}
    \int_{\Omega} (\partial_t  S_{\eps,N}) e_l dx + \int_{\partial_N \Omega} \mathcal{R} e_l d\Gamma 
   & - \int_{\Omega} a_\eps(S_{\eps,N}) P_{c,\eps}'(S_{\eps,N}) \nabla S_{\eps,N} \cdot \nabla e_l dx \\
     & + \int_{\Omega} a_\eps(S_{\eps,N}) \nabla \partial_t \beta_{\eps}(S_{\eps,N}) \cdot \nabla e_l dx = 0,
  \end{split}
  \label{eq-16}
 \end{equation}
for all $l \in \{1,\ldots,N\}$ satisfying the initial condition
\begin{align}
\beta_\eps(S_{\eps,N})(0) = \beta_\eps(S_D)_{\vert_{t=0}} 
    + \Pi_N \big(\beta_\eps(S_i) - \beta_\eps(S_D)_{\vert_{t=0}}\big),
\label{eq-17}                                                                 
\end{align}
where $\Pi_N$ is the orthogonal projector in $L^2(\Omega)$ on the space $V_N$.

\begin{lemma}
 There exists $T_N > 0$ such that problem \eqref{eq-16} has a unique solution which 
 satisfies $ \beta_\eps(S_{\eps,N}) - \beta_\eps(S_D) \in C^1([0,T_N];V_N)$.
\end{lemma}

\noindent{\bf Proof.}
We define the mapping $\Phi_\eps = \beta_\eps^{-1}$ such that 
 $S_{\eps,N} = \Phi_\eps\Big( \beta_\eps(S_D) + \sum_{i=1}^N \gamma_i(t) e_i \Big)$. 
 For simplicity of notation, we omit indices $\eps$ and $N$ from $S_{\eps,N}$.
We can rewrite \eqref{eq-16} as:
\begin{align}
A(\gamma)\frac{d\gamma}{dt} & = B(\gamma)\gamma + F(\gamma, t),\label{eqms-1}\\
\gamma_k(0) &= (e_k,\beta_\eps(S_i)-\beta_\eps(S_D)_{\vert_{t=0}} ), \quad k=1,\dots, N,
\end{align}
where $A \in \Rb^{N\times N}$, $B \in \Rb^{N\times N}$ and $F \in \Rb^N$.
 For the coefficients in \eqref{eqms-1} we have
 \begin{align*}
 A_{l,i} = \int_\Omega a_\eps(S) \nabla e_i \cdot \nabla e_l dx + \int_\Omega \frac{1}{\beta_\eps'(S)} e_i e_l dx, \quad
 B_{l,i}  = \int_\Omega a_\eps(S) P_{c,\eps}'(S) \frac{1}{\beta_\eps'(S)} \nabla e_i \cdot \nabla e_l dx, 
 \end{align*}
 \begin{multline*}
 F_l = -\int_\Omega a_\eps(S) \nabla(\partial_t \beta_\eps(S_D) ) \cdot \nabla e_l dx 
            - \int_\Omega \frac{\beta_\eps'(S_D)}{\beta_\eps'(S)} (\partial_t S_D) e_l dx \\
    + \int_\Omega a_\eps(S)P_{c,\eps}'(S) \frac{\beta_\eps'(S_D)}{\beta_\eps'(S)} \nabla S_D \cdot \nabla e_l dx
    - \int_{\partial_N \Omega} \mathcal{R} e_l d\Gamma .
 \end{multline*}
 
 The matrix $A = A(\gamma)$ is symmetric, positive definite, depending smoothly on $\gamma$. Namely, one has
 \begin{align*}
  A\xi \cdot \xi = \sum_{i,l=1}^N A_{l,i} \xi_i \xi_l = \int_\Omega a_\eps(S)|\nabla u_\xi|^2 dx 
  + \int_\Omega \frac{1}{\beta_\eps'(S)} |u_\xi|^2 dx > 0,
 \end{align*}
where $u_\xi = \sum_i \xi_i e_i$ and  regularized coefficient satisfies
$a_\eps(S) \geq m_a^\eps > 0$ and $\beta_\eps'(S) = \tau_\eps(S)\geq m_\tau^\eps > 0$.
Furthermore, since $F = F(\gamma,t)$ and $B=B(\gamma)$ are continuously differentiable functions of $\gamma$
and continuous functions of $t$, the Cauchy-Lipschitz theorem 
gives the local existence of the solution. $\hfill \Box$

\vspace{2mm}

 {\sf Step 2: a priori estimates uniform with respect to $N$ for the regularized finite-dimensional problem.}

\begin{proposition}
 Under assumptions {\bf (H1)}--{\bf (H5)}, the solution of the approximate problem
 \eqref{eq-16}--\eqref{eq-17} exists for all times $T < +\infty$ and  $\beta_\eps(S_{\eps,N}) - \beta_\eps(S_D) \in C^1\big( [0,T];V_N \big)$.
 \label{Prop-Galerkin}
\end{proposition}

\noindent{\bf Proof.} 
 Let us show boundedness of the solution, uniform with respect to $N$.
 For the test-function in  \eqref{eq-16}  we take $\partial_t \big( \beta_{\eps}(S_{\eps,N}) - \beta_{\eps}(S_D) \big)$.
 In this way we obtain
 \begin{align*}
      \int_{\Omega} & \frac{1}{\tau_{\eps}(S_{\eps,N})} \Big( \partial_t \beta_{\eps}(S_{\eps,N}) \Big)^2 dx 
    + \int_\Omega a_{\eps}(S_{\eps,N}) | \nabla \partial_t \beta_{\eps}(S_{\eps,N}) |^2 dx \\
   & =
   \int_{\Omega}  \frac{1}{\tau_{\eps}(S_{\eps,N})} \partial_t \beta_{\eps}(S_{\eps,N}) \partial_t \beta_{\eps}(S_{D}) dx  
   + \int_{\partial_N \Omega} \mathcal{R} (\partial_t \beta_{\eps}(S_{D})- \partial_t \beta_{\eps}(S_{\eps,N}) )  d\Gamma \\
   & + \int_{\Omega} a_\eps(S_{\eps,N}) P_{c,\eps}'(S_{\eps,N}) \nabla S_{\eps,N} \cdot \nabla \partial_t \beta_\eps(S_{\eps,N}) dx
   - \int_{\Omega} a_\eps(S_{\eps,N}) P_{c,\eps}'(S_{\eps,N}) \nabla S_{\eps,N} \cdot \nabla \partial_t \beta_\eps(S_D) dx \\
   & + \int_{\Omega} a_{\eps}(S_{\eps,N}) \nabla \partial_t \beta_\eps(S_{\eps,N}) \cdot \nabla \partial_t \beta_\eps(S_D) dx.
  \end{align*}
Using the assumptions on the boundedness of the regularized coefficients \eqref{eq-14}--\eqref{eq-15}, the Cauchy inequality,
the fact that all integrals which contain $S_D$ are bounded by some constant, and the trace theorem,
we get the following estimate
 \begin{multline}
  \frac{1}{2\|\tau\|_\infty}  \int_{\Omega} (\partial_t \beta_\eps(S_{\eps,N}))^2 dx  + \frac{m_a^\eps}{4} 
  \int_{\Omega} |\nabla \partial_t \beta_\eps(S_{\eps,N})|^2 dx \\
  \leq C_\eps + \frac{2 \|aP_c'\|_{\infty}^2}{m_a^\eps} \int_{\Omega} |\nabla S_{\eps,N}|^2 dx
  \leq C_\eps + \frac{2 \|aP_c'\|_{\infty}^2}{m_a^\eps (m_\tau^\eps)^2} \int_{\Omega} |\nabla \beta_\eps(S_{\eps,N})|^2 dx,
\label{eq-21}
\end{multline}
where $C_\eps$ is independent on $N$.
Next, from the expression
\[ \nabla \beta_\eps(S_{\eps,N})(t) -  \nabla \beta_\eps(S_{\eps,N})(0) = \int_0^t \partial_\xi \nabla \beta_\eps(S_{\eps,N}(\xi)) d \xi, \]
by integrating over the set $\Omega$, one gets:
\begin{align}
 \| \nabla \beta_\eps(S_{\eps,N}(t))\|_{L^2(\Omega)}^2 \leq C + 
 2t\int_0^t \| \partial _\xi \nabla \beta_\eps(S_{\eps,N}) \|_{L^2(\Omega)}^2 d\xi,
\label{eq-20}      
\end{align}
where $C_\eps$ is independent of $N$.
From \eqref{eq-21} and \eqref{eq-20} we get 
\begin{multline}
  \frac{1}{2\|\tau\|_\infty}  \int_{\Omega} (\partial_t \beta_\eps(S_{\eps,N}))^2 dx  + \frac{m_a^\eps}{4} 
  \int_{\Omega} |\nabla \partial_t \beta_\eps(S_{\eps,N})|^2 dx \\
  \leq C_\eps + \frac{4 \|aP_c'\|_{\infty}^2}{m_a^\eps \big( m_\tau^\eps \big)^2}\, t \int_0^t \|  
  \nabla \partial _\xi \beta_\eps(S_{\eps,N}) \|_{L^2(\Omega)}^2 d\xi.
  \label{BfGronw}
\end{multline}

Now, the Gronwall inequality gives
\begin{align*}
 \max_{t \in [0,T]} \int_{\Omega} |\nabla (\partial_t \beta_\eps(S_{\eps,N}) |^2 dx dt \leq C_\eps,
\end{align*}
where $C_\eps$ is independent of $N$. Finally, from \eqref{BfGronw} we get
\begin{align}
 \frac{1}{2\|\tau\|_\infty}  \max_{t \in [0,T] }  \int_{\Omega} (\partial_t \beta_\eps(S_{\eps,N}))^2 dx  + \frac{m_a}{4} 
  \max_{t \in [0,T] }  \int_{\Omega} |\nabla \partial_t \beta_\eps(S_{\eps,N})|^2 dx 
  \leq C_\eps.	
\label{eq-23}
\end{align}
  We conclude that, since $\partial_t \beta_\eps(S_{\eps,N})$ is bounded in $L^\infty(0,T; H^1(\Omega))$, uniformly with respect to $N$,
the solution for the approximate problem \eqref{eq-16}-\eqref{eq-17} exists for all
times $T < +\infty$ and one has $\beta_\eps(S_{\eps,N}) - \beta_\eps(S_D) \in C^1([0,T]; V_N)$. \hfill $\Box$

\vspace{2mm}

 From now on we consider problem \eqref{eq-16}--\eqref{eq-17} on $[0,T]$, for fixed $T>0$. The calculations done in the proof of Proposition \ref{Prop-Galerkin}
give the following a priori estimates uniform with respect to $N$.

\begin{corollary}
	\label{coroll:bound}
Under hypotheses {\bf (H1)}--{\bf (H5)} there is a constant $C_\eps$ independent of $N$ such that
 \begin{align}
     \| \nabla S_{\eps,N} \|_{L^\infty(0,T;L^2(\Omega))} & \leq C_\eps, \label{eq-24.0}\\ 
     \| \partial_t S_{\eps,N} \|_{L^\infty(0,T;L^2(\Omega))} & \leq C_\eps, \label{eq-25}\\
    \| \nabla \beta_\eps(S_{\eps,N}) \|_{L^\infty(0,T;L^2(\Omega))} & \leq C_\eps, \label{eq-24}\\ 
    \| \nabla (\partial_t \beta_\eps(S_{\eps,N}) ) \|_{L^\infty(0,T;(L^2)(\Omega)^n)} & \leq C_\eps, \label{eq-26}\\
    \| a_\eps(S_{\eps,N}) \nabla (\partial_t \beta_\eps(S_{\eps,N}) ) \|_{L^\infty(0,T;(L^2)(\Omega)^n)} & \leq C_\eps. \label{eq-26-a}
\end{align}
\end{corollary}
\noindent {\bf Proof}.  Estimate \eqref{eq-26} follows directly from  \eqref{eq-23},
while estimate  \eqref{eq-24} follows from \eqref{eq-20} and \eqref{eq-26}.
Moreover, estimate \eqref{eq-23} gives that
\[ \| \partial_t \beta_\eps(S_{\eps,N}) \|_{L^\infty(0,T;L^2(\Omega))} \leq C,    \]
and, since $\tau_\eps$ is bounded (see \eqref{eq-15}), estimate \eqref{eq-25} follows immediately.
Note that \eqref{eq-26-a} follows from \eqref{eq-26} due to the boundedness of the coefficient $a_\eps$.
The estimate  \eqref{eq-24.0} follows from  \eqref{eq-24} and boundedness  of regularized coefficients given by \eqref{eq-15}. $\hfill \Box$

\vspace{2mm}

{\sf  Step 3: passing to the limit when $N \to \infty$.}

Using the a priori estimates obtained in Corollary~\ref{coroll:bound} we are able to pass to the limit
in the weak formulation when $N \to \infty$. 
 More precisely, under hypotheses {\bf (H1)}--{\bf (H5)}, the solution $S_{\eps,N}$ of the approximate problem
 \eqref{eq-16}-\eqref{eq-17} converges to a function $S_\eps \in H^1(Q_T)$ when $N \to +\infty$ satisfying equation \eqref{WFReg}. 
 Moreover, we have $\partial_t \beta_\eps(S_\eps) \in L^2(Q_T)$.
 Now, it is straightforward to prove that the limit function $S_\eps$ satisfies the weak formulation
 \eqref{WFReg} as well as the initial and boundary conditions. This completes the proof of Theorem~\ref{Thm-Reg}. 
 $\hfill \Box$
 
\section{A priori estimates uniform with respect to $\eps$}\label{Sect.6}

 For clarity and better readability the proof of Theorem \ref{Thm-Main} is written in few steps.
 In this section we state a priori estimates uniform with respect to the regularizing
 parameter $\eps$ needed for passing to the limit when $\eps \to 0$. The needed a priori estimates
 are obtained using two different test-functions in weak formulation \eqref{WFReg}. 
 Furthermore, we prove that the integral of the entropy functional given by expression \eqref{Czisar-Kullback}
 is bounded. 
 
 \vspace{2mm}
 
The key-result for proving the global existence of a weak solution for the regularized problem are the a priori estimates given by the following theorem.
\begin{theorem} {\sf A priori estimates uniform with respect to $\eps$.}   \label{AprEstThm}\\
  Under the hypotheses {\bf (H1)}--{\bf (H6)} there is a constant $C>0$, independent of $\eps$, such that every solution
 $S_\eps$ for the regularized problem \eqref{WFReg}
 satisfies
 \begin{align}
 \| \partial_t \beta_\eps(S_\eps) \|_{L^2(Q_T)} & \leq C, \label{AE.7}\\
  \| \sqrt{a_\eps(S_\eps)} \nabla \partial_t \beta_\eps(S_\eps) \|_{L^2(Q_T)} & \leq C,  \label{AE.2}\\
  \| \nabla \beta_\eps(S_\eps) \|_{L^\infty(0,T;L^2(\Omega))} & \leq C, \label{AE.3}\\
   \big\| \sqrt{\tau_\eps(S_\eps)} \nabla \int_0^{S_\eps} \sqrt{-P_{c,\eps}'(\xi)}d\xi \big\|_{L^2(Q_T)} & \leq C. \label{AE.4}\\
 \sup_t \int_{\Omega} \mathcal{E}_\eps(S_\eps) dx dt & \leq C \label{AE.Entr},
\end{align}
where $\mathcal{E}_\eps$ is given below by expression \eqref{Czisar-Kullback_Eps}.
\end{theorem}
\noindent{\bf Proof.}
We test variational equation \eqref{WFReg}
by using two different test-functions. 

{\sf Test-function 1.}
The first choice for the test-function is given by
\begin{align}
  \varphi = \varphi_\eps(S) = \int_{S_D}^S \frac{\tau_\eps(\xi)}{a_\eps(\xi)} d\xi.
 \label{Test5}
\end{align}
We note that this selection for $\varphi_\eps$ is connected to the regularized entropy, corresponding
to $\big( k_{\eps}(S_\eps) \big)^{-1}$ in \cite{M10}.
From Section \ref{Intro}, Subsection {\sf Key ideas}, we recall that
$\mathcal{E}_\eps''(S) =  \tau_\eps(S)/a_\eps(S)$
where
\begin{align}
 \mathcal{E}_\eps(S) = S \int_{S_D}^S \frac{\tau_\eps(\xi)}{a_\eps(\xi)} d\xi - \int_{S_D}^S \xi \frac{\tau_\eps(\xi)}{a_\eps(\xi)} d\xi,
 \label{Czisar-Kullback_Eps}
\end{align}
so we have $\varphi_\eps(S) = \mathcal{E}_\eps'(S)$.

 Using the test-function given by \eqref{Test5}, from \eqref{WFReg} we obtain the following expression
\begin{equation}
 \begin{split}
   & \int_\Omega \mathcal{E}_\eps(S_\eps(t)) dx  + \frac{1}{2} \int_\Omega |\nabla \beta_\eps(S_\eps)|^2 dx \\
   & + \int_0^t\int_\Omega \Big(-P_{c,\eps}'(S_\eps) \Big) \tau_\eps(S_\eps) |\nabla S_\eps|^2 dx dt 
  +  \int_0^T\int_{\partial_N \Omega} \mathcal{R} \varphi d\Gamma dt \\
    & = \int_\Omega \mathcal{E}_\eps(S_i) dx + \frac{1}{2} \int_\Omega |\nabla \beta_\eps(S_i)|^2 dx 
    - \int_0^t\int_\Omega (S_\eps - S_D) \frac{\tau_\eps(S_D)}{a_\eps(S_D)} \partial_t S_D dx dt \\
    & + \int_0^t\int_\Omega \frac{a_\eps(S_\eps)}{a_\eps(S_D)} (-P_{c,\eps}'(S_\eps)) \tau_\eps(S_D) \nabla S_\eps \cdot \nabla S_D dx dt \\
    & + \int_0^t\int_\Omega \frac{a_\eps(S_\eps)}{a_\eps(S_D)} \nabla \partial_t \beta_\eps(S_\eps) \cdot \nabla \beta_\eps(S_D) dx dt.
 \end{split}
\label{AprEstMid.5}
\end{equation}

We denote
\begin{align*}
 J_1 = \int_0^t\int_\Omega (S_\eps - S_D) \frac{\tau_\eps(S_D)}{a_\eps(S_D)} \big( \partial_t S_D \big) dx dt.
\end{align*}
The Cauchy inequality gives
\begin{align}
    |J_1| \leq \delta_1 \int_0^t\int_\Omega \tau_\eps(S_D)^2(S_\eps - S_D)^2 dx dt 
  + \frac{1}{4\delta_1} \int_0^t\int_\Omega \frac{(\partial_t S_D)^2}{a_\eps(S_D)^2}dx dt.
 \label{IntMid.2}
\end{align}

Using the mean value theorem and the assumption that $\tau_\eps(s)$ is monotone increasing function, it is easy to derive the following
inequality
\begin{align}
  \beta_\eps(S_\eps) - \beta_\eps(S_D) \geq \tau_\eps(S_D) ( S_\eps - S_D). 
  \label{Ineq.1}
\end{align}

After inserting \eqref{Ineq.1} into \eqref{IntMid.2}, using the Poincar\'e inequality and the assumptions that $S_D \geq S_{D_{min}} > 0$
one gets:
\begin{align*}
  J_1 \leq \delta_1 C_p \int_0^t\int_\Omega |\nabla \beta_\eps(S_\eps) - \nabla \beta_\eps(S_D)|^2 dx dt + \frac{C}{4\delta_1},
\end{align*}
where the constant $C_p$ comes from the Poincar\'e inequality.
By rewriting the term
\begin{align*}
  \int_0^t\int_\Omega |\nabla \beta_\eps(S_\eps) - \nabla \beta_\eps(S_D)|^2 dx dt 
      = 2 \int_0^t\int_\Omega |\nabla \beta_\eps(S_\eps)|^2 dx dt
      + 2 \int_0^t \int_\Omega |\nabla \beta_\eps(S_D)|^2 dx dt,
\end{align*}
after using the Cauchy-Schwarz and then the Cauchy inequality we get
 \[ |J_1| \leq 2 C_p \delta_1  \int_0^t \int_\Omega |\nabla \beta_\eps(S_\eps)|^2 dx dt \\
      + 2 C_p \delta_1 \int_0^t \int_\Omega |\nabla \beta_\eps(S_D)|^2 dxdt
     + \frac{C}{4\delta_1}. \]
The third term on the left-hand side in \eqref{AprEstMid.5} 
\[ \int_0^t\int_\Omega \Big(-P_{c,\eps}'(S_\eps) \Big) \tau_\eps(S_\eps)  |\nabla S_\eps|^2 dx dt  \]
we rewrite as
\[ \int_0^t\int_\Omega  \tau_\eps(S_\eps) \Big| \nabla \int_0^{S_\eps} \sqrt{-P_{c,\eps}'(\xi)} d\xi\Big|^2 dx dt. \]

Next, we estimate the integral
\[ |J_2| = \int_0^t\int_\Omega \frac{a_\eps(S_\eps)}{a_\eps(S_D)} \Big(-P_{c,\eps}'(S_\eps)\Big) \tau_\eps(S_D) \nabla S_\eps \cdot \nabla S_D dx dt.  \]
After applying the Cauchy inequality we have
\begin{align*}
  J_2 \leq \delta_2 \|a\|_\infty 
     \int_0^t\int_\Omega \tau_\eps(S_\eps) \big| \nabla \int_0^{S_\eps} \sqrt{-P_{c,\eps}'(\xi)} d\xi \big|^2 dx dt 
       + \frac{1}{4\delta_2} \Big\|\frac{aP_C'}{\tau} \Big \|_\infty \int_0^t\int_\Omega \frac{|\nabla \beta_\eps(S_D)|^2}{a_\eps(S_D)^2}dx dt.  
\end{align*}
 It remains to estimate the integral
\[ J_3 = \int_0^t\int_\Omega \frac{a_\eps(S_\eps)}{a_\eps(S_D)} \nabla \partial_t \beta_\eps(S_\eps) \cdot \nabla \beta_\eps(S_D) dx dt.  \]
Similarly like before, the Cauchy inequality gives
\begin{align*}
 J_3 \leq \delta_3 \int_0^t\int_\Omega a_\eps(S_\eps) |\nabla \partial_t \beta_\eps(S_\eps)|^2 dx dt 
    + \frac{\|a\|_\infty}{4\delta_3} \int_0^t\int_\Omega \frac{|\nabla \beta_\eps(S_D)|^2}{a_\eps(S_D)^2}dx dt.
\end{align*}
Finally, concerning the estimate of the integral in \eqref{AprEstMid.5} which contains the boundary-flux term 
$\mathcal{R}$, due to assumption ${\bf (H5)}$, we firstly note that
\begin{align*}
 \Big| \mathcal{R} \int_{S_D}^{S_\eps} 
 \frac{\tau_\eps(\xi)}{a_\eps(\xi)} d\xi \Big| & = \Big| R_0 \, \sigma(S_\eps) \int_{S_D}^{S_\eps} \frac{\tau_\eps(\xi)}{a_\eps(\xi)} d\xi \Big| \\
  & \leq |R_0|  \max_{0 \leq S \leq 1} \sigma(S) \int_A \frac{\tau_\eps(\xi)}{a_\eps(\xi)} d\xi 
   \leq |R_0|  \max_{0 \leq S \leq 1} \sigma(S) \max_{S \in A} \frac{\tau_\eps(S)}{a_\eps(S)} \leq C,
\end{align*}
where we denoted by $\supp$ the support of the function and we take the set $A \subset (0,1)$ such that $\supp S_D \cup \supp(\sigma) \subseteq A$.
In this way we conclude that the boundary-flux integral
\begin{align*}
 \int_0^T\int_{\partial_N \Omega} \mathcal{R} \int_{S_D}^S \frac{\tau_\eps(\xi)}{a_\eps(\xi)} d\xi d\Gamma dt \leq C
\end{align*}
 is bounded by a positive constant independent of $\eps$.

\vspace{2mm}

In this way, from \eqref{AprEstMid.5} we obtain:
\begin{align*}
  \int_\Omega \mathcal{E}_\eps(S_\eps(t)) dx  & + \Big( \frac{1}{2} - 2 C_p\delta_1 \Big) \int_\Omega |\nabla \beta_\eps(S_\eps)|^2 dx \\
      & + \Big( 1 - \delta_2 \| a \|_\infty  \Big)
       \int_0^t\int_\Omega \tau_\eps(S_\eps) |\nabla \int_0^{S_\eps} \sqrt{-P_{c,\eps}'(\xi)} d\xi|^2 dx dt \\
    & \leq \int_\Omega \mathcal{E}_\eps(S_i(t)) dx + \frac{1}{2} \int_\Omega |\nabla \beta_\eps(S_i)|^2 dx  
    + \delta_3 \int_0^t\int_\Omega a_\eps(S_\eps) |\nabla \partial_t \beta_\eps(S_\eps)|^2 dx dt \\
    & + C(\delta_1,\delta_2,\delta_3),
\end{align*}
where $C(\delta_1,\delta_2,\delta_3)$ is independent of $\eps$.
 We note that, using hypothesis {\bf (H6)} and explicit calculation of $\mathcal{E}_\eps(S)$ given in Appendix, we easily
 see that
 \[  \int_\Omega \mathcal{E}_\eps(S_i(x)) dx \leq C, \]
 with $C$ independent of $\eps$.
After taking the supremum with the respect to $t \in [0,T]$ and by choosing appropriate $\delta_i$, we get
\begin{align}
  \sup_{t \in [0,T]} \int_\Omega \mathcal{E}_\eps(S_\eps(t)) dx & + \sup_{t \in [0,T]} \int_\Omega |\nabla \beta_\eps(S_\eps(t))|^2 dx 
    + \int_0^T \int_\Omega \tau_\eps(S_\eps) |\nabla \int_0^{S_\eps} \sqrt{-P_{c,\eps}'(\xi)} d\xi|^2 dx dt \nonumber\\
    & \leq C(\delta) + \delta \int_0^T \int_\Omega a_\eps(S_\eps) |\nabla \partial_t \beta_\eps(S_\eps)|^2 dx dt.
   \label{AprEstMid.5.2.0}
\end{align}

\vspace{2mm}

{\sf Test-function 2.}
 Another choice for the test-function is given by
\begin{align}
  \varphi_\eps(S_\eps) = \partial_t \beta_\eps(S_\eps) - \tau_\eps(S_\eps) \partial_t S_D.
 \label{Test2}
\end{align}
 In this way, from \eqref{WFReg} we get:
\begin{equation}
 \begin{split}
  \int_0^T\int_\Omega \tau_\eps(S_\eps) (\partial_t S_\eps)^2 dx dt & + \int_0^T\int_\Omega a_\eps(S_\eps) |\nabla \partial_t \beta_\eps(S_\eps)|^2 dx dt
  \\ & = \int_0^T \int_\Omega \tau_\eps(S_\eps) (\partial_t S_\eps) (\partial_t S_D) dx dt  \\
   & + \int_0^T \int_\Omega a_\eps(S_\eps) P_{c,\eps}'(S_\eps) \nabla S_\eps \cdot \nabla \partial_t \beta_\eps(S_\eps) dx dt \\
    & - \int_0^T \int_\Omega a_\eps(S_\eps) P_{c,\eps}'(S_\eps) \nabla S_\eps \cdot \nabla \Big( \tau_\eps(S_\eps) \partial_t S_D \Big) dx dt \\
    & + \int_0^T \int_\Omega a_\eps(S_\eps) \nabla \partial_t \beta_\eps(S_\eps) \cdot \nabla \Big( \tau_\eps(S_\eps) \partial_t S_D \Big) dx dt - I_5,
 \end{split}
\label{MidStep.2}
\end{equation}
where $\di I_5 = \int_0^T\int_{\partial_N \Omega} \mathcal{R} \big( \partial_t \beta_\eps(S_\eps) - \tau_\eps(S_\eps) \partial_t S_D \big) d\Gamma dt$.

 We start by estimating the boundary-flux term given by $I_5$. We note that the usage of the trace theorem would give the unpleasant
term $\di \int_0^T\int_\Omega |\nabla \partial_t \beta_\eps(S_\eps)|^2 dx dt$ which would make troubles.
For that reason we use assumption {\bf (H5)}. Firstly, we note that
\begin{align}
 R_0(x,t) \sigma(S_\eps) \partial_t \beta_\eps(S_\eps) = \frac{d}{dt} \Big( R_0 \int_0^{S_\eps} \sigma(\xi) \tau_\eps(\xi) d\xi \Big) 
    - \frac{\partial R_0}{\partial t} \int_0^{S_\eps} \sigma(\xi) \tau_\eps(\xi) d\xi.
\end{align}
Now, we get
 \begin{align*}
  & \Big| \int_0^T\int_{\partial_N \Omega} \mathcal{R} \big( \partial_t \beta_\eps(S_\eps) - \tau_\eps(S_\eps) \partial_t S_D  \big) d\Gamma dt \Big|  \\
  & \leq  \Big| \int_0^T\int_{\partial_N \Omega} R_0(x,t) \sigma(S_\eps) \partial_t \beta_\eps(S_\eps) d\Gamma dt  \Big|
       +  \Big| \int_0^T \int_{\partial_N \Omega} R_0(x,t) \sigma(S_\eps) \tau_\eps(S_\eps) \partial_t S_D d\Gamma dt  \Big|\\
      & \leq  \int_{\partial_N \Omega} |R_0(x,T)| \int_0^1 |\sigma(\xi)| \tau_\eps(\xi) d\xi 
            +  \int_{\partial_N \Omega} | R_0(x,0) | \int_0^1 |\sigma(\xi)| \tau_\eps(\xi) d\xi \\
          &  + \int_0^T \int_{\partial_N \Omega} \Big|\frac{\partial R_0}{\partial t} \Big| \int_0^1 |\sigma(\xi)| \tau_\eps(\xi) d\xi  
            + \|\tau\|_\infty  \int_0^T\int_{\partial_N \Omega} | R_0(x,t)\sigma(S_\eps) | \partial_t S_D d\Gamma dt\leq C.
\end{align*}
  Other integrals on the right-hand side in \eqref{MidStep.2} can be estimated standardly, using the boundedness assumptions 
  of the regularized coefficients and the Cauchy inequality.  In this way we obtain the following expressions:
\begin{multline*}
 \int_0^T\int_\Omega \tau_\eps(S_\eps) (\partial_t S_\eps) (\partial_t S_D) dx dt 
   \leq \delta_1 \int_0^T\int_\Omega \tau_\eps(S_\eps) (\partial_t S_\eps)^2 dx dt 
         + \frac{1}{4\delta_1} \int_0^T\int_\Omega \tau_\eps(S_\eps)  (\partial_t S_D)^2 dx dt,
\end{multline*}
\begin{multline*}
 \int_0^T \int_\Omega a_\eps(S_\eps)  P_{c,\eps}'(S_\eps) \nabla S_\eps \cdot \nabla  \partial_t \beta_\eps(S_\eps) dx dt \\
 \leq  \frac{1}{4\delta_2} \Big\|\frac{aP_c'}{\tau} \Big\|_\infty \int_0^T\int_\Omega \tau_\eps(S_\eps)|\nabla \int_0^{S_\eps} 
  \sqrt{-P_{c,\eps}'(\xi)} d\xi |^2 dx dt \\
      + \delta_2 \int_0^T \int_\Omega a_\eps(S_\eps) |\nabla \partial_t \beta_\eps(S_\eps)|^2 dx dt,
\end{multline*}
\begin{multline*}
 \int_0^T \int_\Omega a_\eps(S_\eps)  \big(-P_{c,\eps}'(S_\eps)\big) \nabla S_\eps \cdot \nabla \Big( \tau_\eps(S_\eps) \partial_t S_D \Big) dx dt \\
    =   \int_0^T \int_\Omega a_\eps(S_\eps) \big(-P_{c,\eps}'(S_\eps)\big)\nabla S_\eps \cdot
   \Big(\tau_\eps(S_\eps) \nabla \partial_t S_D + (\partial_t S_D) \nabla \tau_\eps(S_\eps) \Big) dx dt \\
    = G_1 + G_2.
\end{multline*}
Using the Cauchy inequality, we get 
\begin{align*}
 G_1 \leq \delta_3 \|a P_c'\|_\infty \int_0^T\int_\Omega |\nabla \beta_\eps(S_\eps)|^2 dx dt +
           \frac{1}{4\delta_3} \int_0^T\int_\Omega |\nabla \partial_t S_D|^2 dx dt,
\end{align*}
and
\begin{align*}
 G_2 \leq \|\partial_t S_D\|_{L^\infty(Q_T)} \Big\| \frac{a\tau'}{\tau} \Big\|_\infty \int_0^T\int_\Omega \tau_\eps(S_\eps) 
     \Big|\nabla \int_0^{S_\eps} \sqrt{-P_{c,\eps}(\xi)} d\xi\Big|^2 dx dt. 
\end{align*}
We continue with estimation of the integrals on the right-hand side in \eqref{MidStep.2}. We have:
\begin{multline*}
  \int_0^T \int_\Omega a_\eps(S_\eps) \nabla \partial_t \beta_\eps(S_\eps) \cdot \nabla \Big( \tau_\eps(S_\eps) \partial_t S_D \Big) dx dt \\
   = \int_0^T \int_\Omega a_\eps(S_\eps) \nabla \partial_t \beta_\eps(S_\eps) \cdot \tau_\eps(S_\eps) \nabla \partial_t S_D dx dt +
      \int_0^T \int_\Omega a_\eps(S_\eps) \nabla \partial_t \beta_\eps(S_\eps) \cdot \nabla S_\eps  (\partial_t S_D) \tau'_\eps(S_\eps)dx dt \\
       = K_1 + K_2.
\end{multline*}
Again, using the Cauchy inequality we get
\begin{align*}
 K_1 \leq \delta_4 \int_0^T \int_\Omega a_\eps(S_\eps) |\nabla \partial_t \beta_\eps(S_\eps)|^2 dx dt + 
       \frac{1}{4 \delta_4} \| a \tau^2\|_\infty \int_0^T \int_\Omega |\nabla \partial_t S_D|^2 dx dt,
\end{align*}
and
\begin{align*}
 K_2 \leq \|\partial_t S_D\|_{L^\infty(Q_T)} \Big\| \frac{\sqrt{a}\tau'}{\tau} \Big\|_\infty \Bigg( \delta_5  
 \int_0^T \int_\Omega a_\eps(S_\eps) |\nabla \partial_t \beta_\eps(S_\eps)|^2 
 + \frac{1}{4\delta_5} \int_0^T \int_\Omega |\nabla \beta_\eps(S_\eps)|^2 dx dt \Bigg).
\end{align*}
 By the boundedness assumption of the regularized coefficient $\tau_\eps$ and
 after choosing appropriate constants $\delta_i$ we get:
\begin{align*}
     \frac{1}{2}\int_0^T\int_{\Omega} & \tau_\eps(S_\eps) (\partial_t S_\eps)^2 dx dt 
  +  \frac{1}{2} \int_0^T\int_\Omega a_\eps(S_\eps)|\nabla \partial_t \beta_\eps(S_\eps)|^2 dx dt \\
 & \leq \frac{1}{2} \|\tau\|_\infty \int_0^T \int_\Omega (\partial_t S_D)^2 dx dt \\
 & + \Bigg( \frac{3}{2}  + \|\partial_t S_D\|_\infty \Big\| \frac{\tau' a}{\tau} \Big\|_\infty   \Bigg) 
  \int_0^T\int_{\Omega} \tau_\eps(S_\eps) \Big|\nabla \int_0^{S_\eps} \sqrt{-P_{c,\eps}'(\xi)} d\xi \Big|^2 dx dt\\
 & + \Bigg( \| a P_c'\|_\infty + \frac{3}{2} \|\partial_t S_D\|^2_\infty \Big\| \frac{\tau' \sqrt{a}}{\tau} \Big\|^2_\infty   \Bigg) 
 \int_0^T \int_\Omega |\nabla \beta_\eps(S_\eps)|^2 dx dt + C,
\end{align*}
where the last constant comes from the estimate for the boundary-flux.
 Finally, we have:
\begin{multline}
     \int_0^T\int_{\Omega} \tau_\eps(S_\eps) (\partial_t S_\eps)^2 dx dt 
 +  \int_0^T\int_\Omega a_\eps(S_\eps)|\nabla \partial_t \beta_\eps(S_\eps)|^2 dx dt \\
     \leq  C \Big( 1 + \int_0^T\int_{\Omega} \tau_\eps(S_\eps) \Big|\nabla \int_0^{S_\eps} \sqrt{-P_{c,\eps}'(\xi)} d\xi \Big|^2 dx dt
  + \int_0^T \int_\Omega |\nabla \beta_\eps(S_\eps)|^2 dx dt \Big).
   \label{AprEst.2.1.0}
\end{multline}
We now combine estimate \eqref{AprEst.2.1.0} with previously obtained estimate \eqref{AprEstMid.5.2.0} to get
\begin{multline*}
  \sup_{t \in [0,T]} \int_\Omega  \mathcal{E}_\eps(S_\eps(t)) dx
     + \sup_{t \in [0,T]} \int_\Omega |\nabla \beta_\eps(S_\eps)|^2 dx 
       + \int_0^T\int_\Omega \tau_\eps(S_\eps) |\nabla \int_0^{S_\eps} \sqrt{-P_{c,\eps}'(\xi)} d\xi|^2 dx dt \\
   \leq C(\delta) + C \delta \Big[ 1 + \int_0^T \int_\Omega \tau_\eps(S_\eps) |\nabla \int_0^{S_\eps} \sqrt{-P_{c,\eps}'(\xi)} d\xi|^2 dx dt
    + \int_0^T \int_\Omega |\nabla \beta_\eps(S_\eps)|^2 dx dt \Big],
\end{multline*}
from where it directly follows that
\begin{align*}
  \sup_{t \in [0,T]} \int_\Omega \mathcal{E}_\eps(S_\eps(t)) dx + \sup_{t \in [0,T]} \int_\Omega |\nabla \beta_\eps(S_\eps)|^2 dx 
  + \int_0^T \int_\Omega \tau_\eps(S_\eps) |\nabla \int_0^{S_\eps} \sqrt{-P_{c,\eps}'(\xi)} d\xi|^2 dx dt \leq C.
\end{align*}
 The last inequality gives estimates \eqref{AE.3}, \eqref{AE.4} and \eqref{AE.Entr}.
Furthermore, from \eqref{AprEst.2.1.0} it follows
\begin{align*}
     \int_0^T\int_{\Omega} \tau_\eps(S_\eps) (\partial_t S_\eps)^2 dx dt 
   +  \int_0^T\int_\Omega a_\eps(S_\eps)|\nabla \partial_t \beta_\eps(S_\eps)|^2 dx dt \leq C.
\end{align*}
In this way we get estimate \eqref{AE.2} and also 
\begin{align}
 \| \sqrt{\tau_\eps(S_\eps)} (\partial_t S_\eps) \|_{L^2(Q_T)} & \leq C.
 \label{AE.1}
\end{align}

Finally, estimate \eqref{AE.7} follows directly from \eqref{AE.1} and the fact that, using assumption \eqref{eq-15}, one has 
\begin{align*}
  \int_0^T\int_{\Omega} & \tau_\eps(S_\eps) (\partial_t S_\eps)^2 dx dt 
   \geq \frac{1}{\|\tau\|_\infty}  \int_0^T\int_{\Omega} \big( \partial_t \beta_\eps(S_\eps)  \big)^2 dx dt.
\end{align*}
In this way we conclude the proof of Theorem~\ref{AprEstThm}. $\hfill \Box$

\section{A priori estimate of the mixed-derivative term}\label{Sect.7}

 The main result of this section is given by Proposition \ref{Prop_ThirdOrder} 
where we estimate the third-order term $ \di \nabla( \partial_t \beta_\eps(S_\eps))$.
For that we need two technical results given in Proposition \ref{Prop_Tech} and Lemma \ref{Technical} stated below.
The needed integrability will be achieved by combining between 
results \eqref{Res.1} and \eqref{Res.2} with \eqref{AE.5} and \eqref{AE.6} 
given by the following proposition.
\begin{proposition} \label{Prop_Tech}
 Under hypotheses {\bf (H1)}--{\bf (H6)} there is a constant $C>0$, independent of $\eps$, and constants $\gamma>2$ and $\lambda > 2$
 such that every solution
 $S_\eps$ of the regularized problem \eqref{WFReg}
 satisfies
 \begin{align}
 \| (Z_\eps(S_\eps))^{2-\gamma} \|_{L^\infty(0,T;L^1(\Omega))} & \leq C,\label{AE.5} \\
 \| \big( 1-Z_\eps(S_\eps)\big)^{2 - \lambda} \|_{L^\infty(0,T;L^1(\Omega))} & \leq C. \label{AE.6} 
\end{align}
\end{proposition}
\noindent{\bf Proof.} The result of this proposition follows directly from estimate \eqref{AE.Entr} and
Lemma~\ref{Technical} given below. 
\begin{lemma} \label{Technical}
  For $\gamma>2$ and $\lambda > 2$ there exist positive constants $C$ and $D$,
 independent of $\eps$, such that for any $s \in \Rb$
 \begin{align}
 \mathcal{E}_\eps(s) \geq \mathcal{E}_\eps^0(s) :=  C \Bigg[ \frac{1}{\big(Z_\eps(s) \big)^{\gamma-2}} 
 + \frac{1}{\big( 1- Z_\eps(s) \big)^{\lambda - 2}}   \Bigg] - D, 
 \label{CK_0}
 \end{align}
\end{lemma}

The proof of Lemma~\ref{Technical} is technical and it is given in Appendix.

\vspace{2mm}

Now we state the result which gives the a priori estimate of the third-order term.
\begin{proposition}  \label{Prop_ThirdOrder}
 Let $n \leq 3$ and let us suppose hypotheses {\bf (H1)}--{\bf (H6)}. Furthermore, let the exponents $\beta_1$, $\gamma$, $\mu$,
and $\lambda$ satisfy:
\begin{align}
   & 5 < \beta_1 \leq \gamma < \mu < \frac{5}{6}\gamma + \frac{1}{2}\big( \beta_1 - \frac{10}{3}  \big),
       \quad 5 < \beta_2 \leq \lambda < 3\beta_2-10 \quad
     \textrm{ for } \; n = 3, \label{Add.1}\\
   & 4 <  \beta_1 \leq \gamma < \mu < \sigma \gamma + \frac{1}{2} \big( \beta_1 - 4 \sigma \big), \quad 4 < \beta_2 \leq \lambda \quad
    \textrm{ for } \; n = 2, \label{Add.2} \\
   &  4 < \beta_1 \leq \gamma < \mu < \gamma + \frac{1}{2} \big( \beta_1 - 4  \big), \quad \; \; \; \; 4 < \beta_2 \leq \lambda \quad
    \textrm{ for } \; n = 1, \label{Add.3}
\end{align}
where $\sigma = (p-1)/p$, for arbitrary $p > 2$.
 Then there exist constant $C$, independent of $\eps$, such that any solution $S_\eps$ of the regularized problem \eqref{WFReg} satisfies
 \begin{align}
  \| \nabla (\partial_t \beta_\eps(S_\eps)) \|_{L^{m_0}(Q_T)} \leq C, \; \textrm{ with } \; m_0 \in (1,2).
  \label{ThirdOrder}
 \end{align}
\end{proposition}

\noindent{\bf Proof.}
We start with the estimate
\begin{equation}
 \begin{split}
  & \int_0^T\int_\Omega  |\nabla \partial_t \beta_\eps(S_\eps)|^m dx dt  = 
     \int_0^T\int_\Omega \Big( a_\eps(S_\eps)^{1/2} |\nabla \partial_t \beta_\eps(S_\eps)| \Big)^m \Big(a_\eps(S_\eps)\Big)^{-m/2} dx dt \\
 & \leq \Big( \int_0^T\int_\Omega  \Big( a_\eps(S_\eps)^{1/2}  | \nabla \partial_t \beta_\eps(S_\eps) |\Big)^{mr} dx dt \Big)^{1/r}  
         \Big( \int_0^T\int_\Omega \big( a_\eps(S_\eps)^{-m/2} \big)^{r'} dx dt  \Big)^{1/r'} \\
 & =  \Big( \int_0^T\int_\Omega a_\eps(S_\eps) |\nabla \partial_t \beta_\eps(S_\eps)|^2 dx dt \Big)^{m/2} 
         \Big( \int_0^T\int_\Omega \big( a_\eps(S_\eps) \big)^{-\frac{m}{2-m}} dx dt  \Big)^{1-\frac{m}{2}},
 \end{split}
\label{ExprThOr}
\end{equation}
for $1 < m < 2$, where we used the H\"older inequality with $r = 2/m$ and $r'=2/(2-m)$.
Due to a priori estimate \eqref{AE.2}, the first integral on the right-hand side of \eqref{ExprThOr} is bounded.
We need to conclude that 
\[a_\eps(S_\eps)^{-\frac{m}{2-m}} \in L^1(Q_T), \]
and uniformly bounded with respect to $\eps$.
 Let us set $\delta = m/(2-m)$ and note that
 $1 < \delta < +\infty$. 
We consider the integral
\[ I_\eps = \int_0^T\int_\Omega \big( a_\eps(S_\eps) \big)^{-\delta} dx dt. \]
Now we recall the regularization of the capillary diffusion coefficient $a$ given by \eqref{eq-13}, i.e.
$\di a_\eps(s) = a(Z_\eps(s)) $
where
$\di Z_{\eps}(s)  = (1-2\eps)Z(s) + \eps$.
Using assumption {\bf (H1)}, we get
\begin{align*}
 I_\eps = \iint_{Q_T}  \Bigg\{ \frac{1}{ \big( Z_\eps(S_\eps)\big)^{\mu}} 
          + \frac{1}{\big( 1- Z_\eps(S_\eps) \big)^{\lambda}} \Bigg\}^{\delta} dx dt. 
\end{align*}
 It follows that it is sufficient to undertake to estimate the integrals:
\[ \iint_{Q_T} \big( Z_{\eps}(S_\eps) \big)^{-\mu \delta} dx dt \quad \textrm{ and } \quad 
                \iint_{Q_T} \big( 1 - Z_{\eps}(S_\eps) \big)^{-\lambda \delta}  dx dt.    \]
 In order to get the estimates of the mentioned integrals, we will use the informations about the 
 capillary pressure given by assumption {\bf (H3)}.
 Moreover, from  {\bf (H3)} we have that $\mu \geq \beta_1$ and $\lambda \geq \beta_2$.
   Next, we use a priori estimate \eqref{AE.4} and we have
\begin{align}
  \iint_{Q_T} \tau(Z_\eps(S_\eps)) | P_{c}'(Z_\eps(S_\eps)) | |\nabla Z_\eps(S_\eps)|^2 dx dt \leq C.
  \label{Third.2}
 \end{align}
  
 For simplicity we denote $ y = Z_\eps(S_\eps)$.
 Using assumptions {\bf (H2)} and {\bf (H3)}, expression \eqref{Third.2} can be written as follows
 \begin{align*}
  \iint_{Q_T} \frac{y^{\mu}}{y^\mu + (1-y)^\lambda} 
    \Big[ \TM + \frac{(1-y)^{\lambda}}{y^{\gamma}} \Big]  \Big( \frac{g(y)}{y^{\beta_1}} + \frac{h(y)}{(1-y)^{\beta_2}} \Big) |\nabla y|^2 dx dt \leq C.
 \end{align*}
Since $\di \Big( y^{\mu} + (1-y)^{\lambda} \Big)^{-1} \geq C$
we get
 \begin{align*}
  \iint_{Q_T} \Big( \TM y^{\mu} + y^{\mu - \gamma} (1-y)^{\lambda}  \Big)
     \Big( \frac{g(y)}{y^{\beta_1}} + \frac{h(y)}{(1-y)^{\beta_2}} \Big)  |\nabla y|^2 dx dt \leq C,
 \end{align*}
 from where, using the boundedness from below of functions $g$ and $h$, we get in particular the following estimates:
 \begin{enumerate}
  \item[(i)] $\di \iint_{Q_T} y^{\mu} (1-y)^{-\beta_2} |\nabla y|^2 dxdt \leq C$,
  \item[(ii)] $\di \iint_{Q_T} y^{\mu - \gamma - \beta_1} (1-y)^{\lambda} |\nabla y|^2 dx dt \leq C$.
 \end{enumerate}
We start with expression (ii). By defining
 $\widehat{y} = \min \{y, 1/2 \}$ from (ii) it follows
\begin{align*}
 \Big( \frac{1}{2} \Big)^{\lambda} \iint_{\widetilde{Q}_T} \widehat{y}^{\mu - \gamma - \beta_1} |\nabla \widehat{y}|^2 dx dt \leq C,
\end{align*} 
 so we get
 \begin{align*}
   \iint_{Q_T} \Big| \nabla \big( \widehat{y}^{\alpha_1} \big) \Big|^2 dx dt \leq C,
\end{align*} 
with $\alpha_1 = 1 + (\mu - \gamma - \beta_1)/2$.
In this way we conclude that
\begin{align*}
 \Big( \widehat{Z_\eps(S_\eps)} \Big)^{\alpha_1}
   = \Big( \min \{ Z_\eps(S_\eps), \frac{1}{2} \} \Big)^{\alpha_1} \in L^2(0,T;H^1(\Omega))
\end{align*}
 is uniformly bounded.
 Because of the Sobolev embedding theorem, one has uniform boundedness in the following spaces:
\begin{align}
   \Big( Z_\eps(S_\eps) \Big)^{\alpha_1} \in 
    \begin{cases}
   L^2(0,T;C(\overline{\Omega})) & \; \textrm{ for } \; n=1, \\
   L^2(0,T;L^p(\Omega)),\, \forall p \in (1,+\infty) & \; \textrm{ for } \; n=2, \\
   L^2(0,T;L^6(\Omega)), & \; \textrm{ for } \; n=3.
 \end{cases}
  \label{Res.1}
 \end{align}
  In similar way, from expression (i) we get uniform boundedness in the following spaces:
 \begin{align}
  \Big( 1 - Z_\eps(S_\eps) \Big)^{\alpha_2} \in
  \begin{cases}
   L^2(0,T;C(\overline{\Omega})) & \; \textrm{ for } \; n=1, \\
   L^2(0,T;L^p(\Omega)),\, \forall p \in (1,+\infty) & \; \textrm{ for } \; n=2, \\
   L^2(0,T;L^6(\Omega)), & \; \textrm{ for } \; n=3,
  \label{Res.2}
 \end{cases}
 \end{align}
 where $\alpha_2 = 1 - \beta_2/2$.

\vspace{2mm}

 In the following part of the proof of Proposition~\ref{Prop_ThirdOrder}, using  \eqref{Res.1} and \eqref{Res.2} we are going to prove that for some $\delta > 1$
 the term $( a_\eps(S_\eps) )^{-\delta}$ is  uniformly bounded in $L^1(Q_T)$.
 For simplicity, we continue to use notation $y = Z_\eps(S_\eps)$. 
 Again, we will consider separately cases $n=1,2$ and $3$.
 
 Let us start with the case $n=3$. We rewrite the expression $\int_\Omega y^{-\mu \delta} dx$ 
 using
 $ -\mu \delta = \alpha_1 \Theta + (2-\gamma) \Theta_1$,
 and then after the H\"older inequality one has
\begin{align*}
 \int_\Omega y^{-\mu \delta} dx & = \int_\Omega y^{\alpha_1 \Theta} \; y^{(2-\gamma)\Theta_1} dx 
    \leq  \Big( \int_\Omega  y^{\alpha_1 \Theta p_1} dx \Big)^{\frac{1}{p_1}} \; \Big( \int_\Omega y^{(2-\gamma)\Theta_1 p_2 } dx \Big)^{\frac{1}{p_2}}.
\end{align*}
We take $p_1 = 6/\Theta$ and $p_2 = 6/(6-\Theta)$. 
Next, for $\Theta = 2$ and $\Theta_1 = 2/3$ we get
\begin{align*}
 \iint_{Q_T} y^{-\mu \delta} dx dt  \leq \int_0^T 
   \Big( \int_{\Omega} y^{6\alpha_1} dx \Big)^{1/3} dt \cdot \max_{0\leq t \leq T} \Big( \int_{\Omega} y^{2-\gamma} dx \Big)^{2/3}.
\end{align*}
Because of \eqref{Res.1}, the first integral on the right-hand side is uniformly bounded. On the other side due to \eqref{AE.5}, 
the second integral is bounded as well. 
From the condition that 
\begin{align}
 \delta = - \frac{1}{\mu} \Big( \frac{10}{3} + \mu - \frac{5}{3} \gamma - \beta_1 \Big) > 1 
 \label{Delta}
\end{align}
we get
\begin{align}
 \frac{10}{3} + 2\mu - \frac{5}{3} \gamma < \beta_1.
 \label{Est.n3.1}
\end{align}
Due to assumption {\bf (H3)}, we have
\begin{align}
 \frac{10}{3} < \beta_1 \leq \gamma < \mu < \frac{5}{6} \gamma + \frac{1}{2} \big( \beta_1 - \frac{10}{3} \big).
 \label{EstCoeff.0}
\end{align}
Using $\mu > \gamma $, from \eqref{Est.n3.1} we have $-10 + 3\beta_1 > \gamma$,
 while $0 < \beta_1 \leq \gamma < \mu$ gives $5\gamma > 3\beta_1 + 10$. By adding two obtained inequalities
 it follows that  $\gamma > 5$.
 Furthermore, since $  \gamma < \mu  $, one has 
$\mu > 5 \; \textrm{ and } \; \beta_1 > 5$.
Finally, for $n=3$ we obtain the following sequence
 of inequalities:
 \begin{align*}
  5 < \beta_1 \leq \gamma < \mu < \frac{5}{6} \gamma + \frac{1}{2} \big( \beta_1 - \frac{10}{3} \big),
 \end{align*}
which gives \eqref{Delta}.

For $n=2$, we write again $-\delta \mu =  \Theta \alpha_1 + (2-\gamma) \Theta_1$.
Next, we take $\Theta = 2$ and $\Theta_1 = 1 - 2/p$ 
and we 
apply the H\"older inequality with $p_1 = p/2$ and $p_2 = p/(p-2)$, where $p>2$ is arbitrary finite.
In this way we get:
\begin{align*}
 \int_0^T \int_{\Omega} y^{-\mu \delta} dx dt & = \int_0^T \int_\Omega y^{2\alpha_1} \; y^{(2-\gamma)\frac{p-2}{p}} dx dt
       \leq \int_0^T \Big( \int_{\Omega} y^{\alpha_1 p} dx \Big)^{2/p} \Big( \int_{\Omega} y^{2-\gamma} dx \Big)^{(p-2)/p} dt\\
      & = \int_0^T \Big( \int_{\Omega} y^{\alpha_1 p} dx \Big)^{2/p} dt \cdot \max_{0\leq t \leq T} \Big( \int_{\Omega} y^{2-\gamma} dx \Big)^{(p-2)/p}.
\end{align*}
Now, taking into account the expression for $\alpha_1$, from $\delta>1$ we get the condition
\[  -\frac{1}{\mu} \Big[ 2 + \mu - \gamma - \beta_1 + (2-\gamma) \frac{p-2}{p}  \Big] > 1. \]
In this way we obtain the following inequality
\begin{align}
 4 \frac{p-1}{p} + 2\mu & - 2\gamma \frac{p-1}{p} <  \beta_1, \label{EstCoeff.1}
 \end{align}
 with $\gamma  > 2$ and $0 < \beta_1  \leq \gamma < \mu$. Next, we denote $ \nu = (p-1)/p$ and from \eqref{EstCoeff.1} we get 
\begin{align}
\mu <  \nu \gamma + \frac{1}{2} \big( \beta_1 - 4\nu \big). 
 \label{EstCoeff.3}
\end{align}
Moreover, the upper bound in \eqref{EstCoeff.3}
is the lowest for $\gamma = \beta_1$, so from
$\beta_1 < \nu \beta_1 + ( \beta_1 - 4\nu )/2$,
we get 
$ \beta_1 > 4\nu/(2\nu-1) =: \tilde{f}(\nu)$.
Since $\tilde{f}(\nu)$ is nonincreasing function for $\nu>1/2$, it is sufficient to have $\beta_1 > \tilde{f}(1) = 4$. 

\vspace{2mm}

For $n=1$, we have:
\begin{align*}
 \int_\Omega y^{-\delta \mu} dx  = 
 \int_\Omega y^{\alpha_1 \Theta} y^{(2-\gamma) \Theta_1} dx \leq \big( \max_\Omega y^{\alpha_1} \big)^2 \int_\Omega y^{2-\gamma} dx
   \in L^1(0,T),
\end{align*}
where we chose $\Theta = 2$ and $\Theta_1 = 1$. The above calculation gives 
$-\delta \mu = 2 \alpha_1 + (2-\gamma)$,
and consequently
$ \delta = -( 4 + \mu - 2\gamma - \beta_1 )/ \mu > 1$,
from where it follows
$ \beta_1 > 4 + 2(\mu - \gamma)$,
which gives $\beta_1 > 4$. 
By taking into account assumption {\bf (H3)}, we conclude that for $n=1$ one has 
\[ 4 < \beta_1 \leq \gamma < \mu < \gamma + \frac{1}{2} (\beta_1 - 4). \]

\vspace{2mm}

 Next, we consider the term $\di (1-Z_\eps(S_\eps))^{-\lambda \delta}$. Similarly like above, 
 using the H\"older inequality and estimates \eqref{AE.6} and \eqref{Res.2}, for $n=1, 2$ and $3$
 we get the expressions for parameters $\beta_2$ and $\lambda$ as it is given by \eqref{Add.1}--\eqref{Add.2}--\eqref{Add.3}.

 Finally, we conclude that $a_\eps^{-m/(2-m)} \in L^1(Q_T)$ for $m=m(n)$. This ends the proof
of Proposition \ref{Prop_ThirdOrder}. $\hfill \Box$

\section{End of proof of Theorem~\ref{Thm-Main}}\label{Proof_Thm_Main}

 For conducting the passage to the limit in \eqref{WFReg} when $\eps \to 0$, we will need almost everywhere convergence
 of the sequence $S_\eps$ towards $S$, as $\eps \to 0$. For obtaining this convergence we need results given by the following proposition
 which, roughly speaking, comes from estimate \eqref{AE.Entr} using the results given in the proof of Lemma~\ref{Technical}.
\begin{proposition}
 Under hypotheses {\bf(H1)}--{\bf(H6)} there exist a constant $C>0$, independent of $\eps$, such that every solution
 $S_\eps$ for the regularized problem \eqref{WFReg}, satisfies:
 \begin{align}
   \| S_\eps^-\|_{L^\infty(0,T;L^2(\Omega))} & \leq C  \eps^{\gamma/2}, \quad |\{ S_\eps \leq 0 \} | \leq C \eps^{\gamma-2}, \label{NegPart}\\
   \| (S_\eps-1)^+)\|_{L^\infty(0,T;L^2(\Omega))} & \leq C \eps^{\lambda/2},\quad  |\{ S_\eps \geq 1 \} | \leq C \eps^{\lambda-2}, \label{PozPart}
 \end{align}
 where $|\cdot|$ represents the measure on $Q_T$.
 \label{Meas}
\end{proposition}

The proof of Proposition~\ref{Meas} is based on the proof of Lemma~\ref{Technical} and it is given in the Appendix, as well. 
 
 \vspace{3mm}
 
From derived a priori estimates \eqref{AE.7}-\eqref{AE.Entr} we conclude that there exists a subsequence of $S_\eps$, denoted by the same
subscript, and $q \in H^1(Q_T)$ 
such that, as $\eps \to 0$ one has:
\[ \beta_\eps(S_\eps) \rightharpoonup q \; \textrm{ weakly in } \; H^1(Q_T). \]
Because of the compact embedding $H^1(Q_T) \hookrightarrow L^2(Q_T)$ it follows
\begin{align}
\beta_\eps(S_\eps) \to q \; \textrm{ strongly in } \; L^2(Q_T)  \; \textrm{ and  a.e. on} \; Q_T
\; \textrm{ on some subsequence. }
 \label{Beta.0}
\end{align}
Let us write $q = q^+ + q^-$. We are going to show that $q^- = 0$.
Using the definition of $\beta$ given in \eqref{beta}, its extension to whole $\Rb$ and the mean-value theorem, direct calculation gives
\begin{align}
 |\beta_\eps(S) - \beta(S)| \leq C \eps \Big( |S^-| + 1 + |(S-1)^+| \Big).
 \label{Beta.1}
\end{align}
From \eqref{Beta.1} and Proposition~\ref{Meas} we get 
\[ \| \beta_\eps(S_\eps) - \beta(S_\eps)\|_{L^2(Q_T)} \leq C \eps. \]
Next, we write $\di \beta_\eps(S_\eps) = \big( \beta_\eps(S_\eps) \big)^+ + \big( \beta_\eps(S_\eps) \big)^-$.
It is obvious that
\[ \| \big( \beta_\eps(S_\eps) \big)^- \|_{L^2(Q_T)} = \tau(\eps)\|S_\eps^-\|_{L^2(Q_T)} \to 0,\; \textrm{ as } \; \eps \to 0, \]
so we conclude that $q \geq 0$ a.e. on $Q_T$.
On the other side, one has 
\[ S_\eps^+ = \beta^{-1}\big((\beta(S_\eps))^+\big) \to \beta^{-1}(q^+) = \beta^{-1}(q)  \; \textrm{ a.e. on }\; Q_T. \]
We define $\di \beta^{-1}(q) =:S$ and we conclude that $S \geq 0$.
In this way we get that
\begin{align}
  S_\eps \to S \; \textrm{ a.e. on } \;  Q_T. 
 \label{S_a.e.}
\end{align}
Using convergence \eqref{S_a.e.} we also have that
$\di (S_\eps - 1)^+ \to (S-1)^+ = 0 \; \textrm{ a.e. on } \;  Q_T$, 
where the last equality comes directly from  \eqref{PozPart}.
In this way we get $S \leq 1$.

Finally, concerning the passage to the limit in \eqref{WFReg}  when $\eps \to 0$ we comment only the integral 
  $\di \int_0^T \int_{\Omega} a_\eps(S_\eps) \nabla \partial_t\beta_\eps(S_\eps) \cdot \nabla \varphi dx dt$.
 From \eqref{AE.2} we conclude that
$\di \sqrt{a_\eps(S_\eps)} \nabla \partial_t \beta_\eps(S_\eps) \rightharpoonup k$ weakly in $L^2(Q_T)$,
 as  $\eps \to 0$.
 Furthermore, from \eqref{ThirdOrder}
 and the fact that $\beta_\eps(S_\eps) \to \beta(S)$ strongly in $L^2(Q_T)$,
 we get 
 $\di \nabla \partial_t \beta_\eps(S_\eps)  \rightharpoonup \nabla \partial_t \beta(S)$
 weakly in $L^{m_0}(Q_T)$, for $m_0 \in (1,2)$ given by Proposition~\ref{Prop_ThirdOrder}. 
 Because of \eqref{S_a.e.}, we get $\sqrt{a_\eps(S_\eps)} \to \sqrt{a(S)} \; \textrm{ a.e.}$
 when $\eps \to 0$. Now, the
Lebesgue theorem on dominated convergence for $L^p$ spaces gives
\begin{align}
 \sqrt{a_\eps(S_\eps)} \to \sqrt{a(S)} \; \textrm{ strongly in }\; L^p(Q_T), \quad \forall p < \infty, \quad \textrm{ as }\; \eps \to 0.
 \label{Sqrt_Lp}
\end{align}
The passage to the limit when $\eps \to 0$ of the integral
$\di \int_0^T \int_{\Omega} a_\eps(S_\eps) \nabla \partial_t\beta_\eps(S_\eps) \cdot \nabla \varphi dx dt$ now 
follows as a product of weak and strong converging sequences. 
The passage to the limit in other terms in \eqref{WFReg} is straightforward.
This ends the proof of Theorem~\ref{Thm-Main}. $\hfill \Box$

\section{Appendix}

{\bf Proof of Lemma~\ref{Technical}.}
 Firstly, we recall the definition of the entropy functional $ \mathcal{E}(s)$ given by expression \eqref{Czisar-Kullback}
 with assumptions on the coefficients $a(s)$ and $\tau(s)$ given by {\bf (H1)} and {\bf (H2)}.
 Let us firstly introduce the notation:
\begin{align*}
 h(s) = \frac{\tau(s)}{a(s)} = h_1^\gamma(s) + \TM h_2^\lambda(s),
\end{align*}
where
$\di h_1^\gamma(s) =  \frac{1}{s^\gamma}$ and  $\di h_2^\lambda(s) = \frac{1}{(1-s)^\lambda}$ with $\gamma > 0, \lambda > 0$.

 In Section~\ref{Sect.4} we regularized the expression for $h(s)$ in the following way:
\begin{align*}
  h_{1,\eps}^\gamma(s)  = \big( Z_\eps(s) \big)^{-\gamma},\quad
  h_{2,\eps}^\lambda(s)  = \big(1-Z_\eps(s)\big)^{-\lambda},
\end{align*}
where $Z_\eps(s)$ was given by
\begin{align*}
 Z_\eps(s) = \left\{
  \begin{array}{ll}
    \eps & \text{if } s < 0,\\
    (1-2\eps)s + \eps & \text{if } 0 \leq s \leq 1,\\
    1-\eps & \text{if } s > 1.
  \end{array} \right.
\end{align*}

Regularized entropy functional reads
\begin{align*}
  \mathcal{E}_\eps(s) =  \mathcal{E}_{1,\eps}^\gamma(s) + \TM \mathcal{E}_{2,\eps}^\lambda(s),
\end{align*}
where
\begin{align*}
  \mathcal{E}_{1,\eps}^\gamma(s) & = s \int_{S_D}^s h_{1,\eps}^\gamma(\xi) d\xi - \int_{S_D}^s \xi h_{1,\eps}^\gamma(\xi) d\xi,\\
  \mathcal{E}_{2,\eps}^\lambda(s) & = s \int_{S_D}^{s} h_{2,\eps}^\lambda(\xi) d\xi - \int_{S_D}^{s} \xi h_{2,\eps}^\lambda(\xi) d\xi.
\end{align*}

 Direct calculation gives the following cases:
 
 {\bf Case 1:}  $0 \leq s \leq 1$
\begin{align*}
 \mathcal{E}_{1,\eps}^\gamma(s) & = \frac{1}{(1-2\eps)^2} \frac{1}{(\gamma-1)(\gamma-2)} \frac{1}{(Z_\eps(s))^{\gamma-2}}
                                   + \frac{1}{(1-2\eps)^2} \frac{1}{(\gamma-1)} \frac{Z_\eps(s)}{(Z_{\eps}(S_D))^{\gamma-1}} \\
                                   & - \frac{1}{(1-2\eps)^2} \frac{1}{\gamma-2}\frac{1}{(Z_{\eps}(S_D))^{\gamma-2}}
                                   \geq \frac{C}{Z_\eps(s)^{\gamma-2}} - D, \\
   \mathcal{E}_{2,\eps}^\lambda(s) & =  \frac{1}{(1-2\eps)^2} \frac{1}{(\lambda-1)(\lambda-2)} \frac{1}{(1-Z_\eps(s))^{\lambda-2}}  
                                   + \frac{1}{(1-2\eps)^2} \frac{1}{(\lambda-1)} \frac{1-Z_\eps(s)}{(1 - Z_{\eps}(S_D))^{\lambda-1}}\\
                                   & - \frac{1}{(1-2\eps)^2} \frac{1}{\lambda-2}\frac{1}{(1 - Z_{\eps}(S_D))^{\lambda-2}}
                                   \geq  \frac{C}{(1-Z_\eps(s))^{\lambda-2}} - D.                                
\end{align*}

 {\bf Case 2:}  $s < 0$
 \begin{align*}
 \mathcal{E}_{1,\eps}^\gamma(s) & =  \frac{s^2}{2 \eps^\gamma} - \frac{1}{(1-2\eps)^2} \frac{1}{\gamma-2} \frac{1}{\big( Z_\eps(S_D) \big)^{\gamma-2}} 
                                    + \frac{s}{1-2\eps} \Big[ \frac{1}{\gamma-1} \frac{1}{\big( Z_\eps(S_D)\big)^{\gamma-1}} - \frac{1}{\eps^{\gamma-1}} \Big] \\
                                   & + \frac{1}{(1-2\eps)^2} \frac{1}{\gamma-1} \frac{\eps}{\big( Z_\eps(S_D)\big)^{\gamma-1}}
                                     + \frac{1}{(1-2\eps)^2} \frac{1}{(\gamma-1)(\gamma-2)} \frac{1}{\eps^{\gamma-2}}. 
\end{align*}
 Using $s/(Z_\eps(S_D))^{\gamma-1} \geq s/\eps^{\gamma-1}$ we can estimate
 \begin{align*}
 \mathcal{E}_{1,\eps}^\gamma(s) & \geq  \frac{s^2}{2 \eps^\gamma} - \frac{1}{(1-2\eps)^2} \frac{1}{\gamma-2} \frac{1}{\big( Z_\eps(S_D) \big)^{\gamma-2}} 
                                    + \frac{|s|}{1-2\eps} \frac{\gamma-2}{\gamma-1} \frac{1}{\eps^{\gamma-1}}  \\
                                   & + \frac{1}{(1-2\eps)^2} \frac{1}{\gamma-1} \frac{\eps}{\big( Z_\eps(S_D)\big)^{\gamma-1}}
                                     + \frac{1}{(1-2\eps)^2} \frac{1}{(\gamma-1)(\gamma-2)} \frac{1}{\eps^{\gamma-2}} \\
                                     & \geq \frac{s^2}{2 \eps^\gamma} +  \frac{|s|}{1-2\eps} \frac{\gamma-2}{\gamma-1} \frac{1}{\eps^{\gamma-1}}
                                        +  \frac{1}{(1-2\eps)^2} \frac{1}{(\gamma-1)(\gamma-2)} \frac{1}{\eps^{\gamma-2}} - D.
 \end{align*}
\begin{align*}
  \mathcal{E}_{2,\eps}^\lambda(s) &  =  \frac{1}{(1-\eps)^\lambda} \frac{s^2}{2} - \frac{1}{(1-2\eps)^2} \frac{1}{\lambda-1} 
                                     \frac{(1-2\eps)s-(1-\eps)}{(1-Z_\eps(S_D))^{\lambda-1}} + \frac{1}{(1-\eps)^{\lambda-1}} \frac{1}{\lambda-1} \frac{s}{(1-2\eps)} \\
                                   & + \frac{s}{1-2\eps} \frac{1}{\lambda-1} \Big[ \frac{1}{(1-\eps)^{\lambda-1}} - \frac{1}{(1-Z_\eps(S_D))^{\lambda-1}} \Big]
                                    + \frac{1}{(1-2\eps)^2} \frac{1}{\lambda-1} \frac{1-\eps}{(1-Z_\eps(S_D))^{\lambda-1}} \\
                                   & + \frac{1}{(1-2\eps)^2} \frac{1}{(1-\eps)^{\lambda-2}} \frac{1}{(\lambda-1)(\lambda-2)}
                                     - \frac{1}{(1-2\eps)^2} \frac{1}{\lambda-2} \frac{1}{(1-Z_\eps(S_D))^{\lambda-2}} \geq -D.                                                        
\end{align*}

 {\bf Case 3:}  $s > 1$
 \begin{align*}
     \mathcal{E}_{1,\eps}^\gamma(s) & = \frac{(1-2\eps)s+\eps}{(1-2\eps)^2} \frac{1}{\gamma-1} 
                       \Big( \frac{1}{Z_\eps(S_D)^{\gamma-1}} - \frac{1}{(1-\eps)^{\gamma-1}} \Big)  \\
                      & - \frac{1}{(1-2\eps)^2} \frac{1}{\gamma-2} \Big( \frac{1}{Z_\eps(S_D)^{\gamma-2}} - \frac{1}{(1-\eps)^{\gamma-2}}  \Big)
                        + \frac{(s-1)^2}{2(1-\eps)^{\gamma}} \geq -D, \\
     \mathcal{E}_{2,\eps}^\lambda(s) & = \frac{(1-2\eps)s + \eps -1}{(1-2\eps)^2} \frac{1}{\lambda -1} 
                                       \Big[ \frac{1}{\eps^{\lambda-1}} - \frac{1}{(1-Z_\eps(S_D))^{\lambda-1}} \Big]  \\ 
                                      & + \frac{1}{(1-2\eps)^2} \frac{1}{\lambda-2} \Big[ \frac{1}{\eps^{\lambda-2}} - \frac{1}{(1-Z_\eps(S_D))^{\lambda-2}} \Big]
                                       + \frac{(s-1)^2}{2\eps^\lambda}.  
\end{align*}
 We note that $ \mathcal{E}_{2,\eps}^\lambda(s)$ can be transformed in the following form:
  \begin{align*}
     \mathcal{E}_{2,\eps}^\lambda(s) & = \frac{(1-2\eps)s + 2\eps -1}{(1-2\eps)^2} \frac{1}{\lambda -1} 
                                       \Big[ \frac{1}{\eps^{\lambda-1}} - \frac{1}{(1-Z_\eps(S_D))^{\lambda-1}} \Big]  \\ 
                                      & + \frac{1}{(1-2\eps)^2} \frac{1}{(\lambda-1)(\lambda-2)} \frac{1}{\eps^{\lambda-2}} 
                                      - \frac{1}{(1-2\eps)^2} \frac{1}{\lambda-2} \frac{1}{(1-Z_\eps(S_D))^{\lambda-2}} \\
                                      & + \frac{\eps}{(1-2\eps)^2} \frac{1}{\lambda-1} \frac{1}{(1-Z_\eps(S_D))^{\lambda-1}}  + \frac{(s-1)^2}{2\eps^\lambda}\\
                                      & \geq \frac{1}{(1-2\eps)^2} \frac{1}{(\lambda-1)(\lambda-2)} \frac{1}{\eps^{\lambda-2}} 
                                      +  \frac{(s-1)^2}{2\eps^\lambda} -D,
 \end{align*}
 where we neglected the first term on the right-hand side since it is positive.
From derived estimates on $\mathcal{E}_{1,\eps}^\gamma(s)$ and $\mathcal{E}_{2,\eps}^\lambda(s)  $  Lemma~\ref{Technical} follows directly. 
$\hfill \Box$
 
  \vspace{3mm}

 {\bf Proof of Proposition~\ref{Meas}.} 
 We start the proof by rewriting the result \eqref{AE.Entr}, i.e. we have the following estimate for a.e. $t \in (0,T)$:
 \begin{align*}
  \int_{\Omega \cap \{S_\eps < 0\}} \mathcal{E}_\eps(S_\eps) dx +  \int_{\Omega \cap \{0 \leq S_\eps \leq 1\}} \mathcal{E}_\eps(S_\eps) dx
   + \int_{\Omega \cap \{S_\eps > 1\}} \mathcal{E}_\eps(S_\eps) dx \leq C.
 \end{align*}
 From previous estimates on $\mathcal{E}_{1,\eps}^\gamma(s)$ and $\mathcal{E}_{2,\eps}^\lambda(s)$ we get
\begin{align*}
 \frac{1}{2\eps^\gamma} \int_{\Omega} (S_\eps^-(t))^2 dx 
 +\frac{1}{(1-2\eps)^2} \frac{1}{(\gamma-1)(\gamma-2)} \frac{1}{\eps^{\gamma-2}} |\{ S_\eps \leq 0 \} | \leq C,
\end{align*}
and 
\[ \frac{1}{2\eps^\lambda} \int_{\Omega} (S_\eps -1)^+ dx 
  + \frac{1}{(1-2\eps)^2} \frac{1}{(\lambda-1)(\lambda-2)} \frac{1}{\eps^{\lambda-2}} |\{ S_\eps \geq 1 \}| \leq C, \]
This completes the proof of Proposition~\ref{Meas}. $\hfill \Box$

\section*{Acknowledgement}
 The author acknowledges support from   
 the Croatian Science Foundation HRZZ-IP-2013-11-3955,
 and the Austrian-Croatian Project of the Austrian Exchange Service (\"OAD) and the Ministry of Science and Education of the Republic of Croatia (MZOS).
 
 This work was partially done during the author was visiting 
 Laboratoire de Math\'e\-ma\-tiques et de leurs Applications at the University of Pau and Pays de l'Adour.
 The author warmly thanks Brahim Amaziane not only for hospitality but also for introducing her to the interesting problems
 concerning the nonequilibrium two-phase flow through porous media. 


\end{document}